\documentclass[12pt]{amsart}
\usepackage{graphicx}
\usepackage{amsfonts}
\usepackage{amssymb}
\usepackage{latexsym}
\usepackage{amscd}

\renewcommand{\marginpar}[1]{}
\catcode`\@=12

\def\Empty{}
\newcommand\oplabel[1]{
  \def\OpArg{#1} \ifx \OpArg\Empty {} \else
  	\label{#1}
  \fi}
		
\long\def\realfig#1#2#3#4{
\begin{figure}[htbp]
\centerline{\includegraphics[width=#4]{#2}}
\caption[#1]{#3}
\oplabel{#1}
\end{figure}}

\newcommand{\comm}[1]{}

\input{psfig}
\newtheorem{thm}{Theorem}[section]
\newtheorem{cor}[thm]{Corollary}
\newtheorem{lem}[thm]{Lemma}
\newtheorem{prop}[thm]{Proposition}

\newtheorem{tet}[thm]{Conditional Implication}
\newtheorem{mlem}[thm]{Main Lemma}

\newenvironment{pf}{\proof[\proofname]}{\endproof}
\newenvironment{pf*}[1]{\proof[#1]}{\endproof}
\usepackage{euscript}

\usepackage[OT2,OT1]{fontenc}

\newcommand{\cal}[1]{{\mathcal #1}}

\newcommand{\beq}{\begin{equation}}
\newcommand{\eeq}{\end{equation}}
\newcommand{\eref}[1]{(\ref{#1})}

\newcommand{\ve}{\varepsilon}
\newcommand{\de}{\delta}
\newcommand{\al}{\alpha}
\newcommand{\be}{\beta}
\newcommand{\ga}{\gamma}

\newcommand{\om}{\omega}

\theoremstyle{definition}
\newtheorem{defn}{Definition}[section]

\theoremstyle{remark}
\newtheorem{rem}{Remark}[section]

\renewcommand{\deg}{\operatorname{deg}}
\newcommand{\riem}{\hat{\CC}}

\newcommand{\dist}{\operatorname{dist}}

\newcommand{\cl}{\operatorname{cl}}

\newcommand{\eps}{\epsilon}

\numberwithin{equation}{section}
\newcommand{\thmref}[1]{Theorem~\ref{#1}}
\newcommand{\propref}[1]{Proposition~\ref{#1}}
\newcommand{\secref}[1]{\S\ref{#1}}
\newcommand{\lemref}[1]{Lemma~\ref{#1}}
\newcommand{\corref}[1]{Corollary~\ref{#1}}

\newcommand{\cM}{{\cal M}}

\newcommand{\cB}{{\cal B}}

\newcommand{\cC}{{\cal C}}
\newcommand{\cH}{{\cal H}}

\newcommand{\cD}{{\cal D}}

\newcommand{\CC}{{\mathbb C}}
\newcommand{\RR}{{\mathbb R}}

\newcommand{\ZZ}{{\mathbb Z}}
\newcommand{\NN}{{\mathbb N}}
\newcommand{\DD}{{\mathbb D}}

\newcommand{\QQ}{{\mathbb Q}}

\begin{document}
\addtolength{\evensidemargin}{-0.7in}
\addtolength{\oddsidemargin}{-0.7in}

\title[Non-computable Julia sets]{Non-computable Julia sets}
\author{M. Braverman, M. Yampolsky}
\thanks{The first author's research is supported by an NSERC CGS scholarship}
\thanks{The second author's research is supported by NSERC operating grant}
\date{\today}
\begin{abstract}
We show that under the definition of computability which is natural
from the point of view of applications, there exist non-computable
quadratic Julia sets.
\end{abstract}
\maketitle

\section{Summary of the paper}

\noindent
Polynomial Julia sets have emerged as the most studied examples of 
fractal sets generated by a dynamical system. Apart from the beautiful
mathematics, one of the reasons for their popularity is the beauty of 
the computer-generated images of such sets. The algorithms used 
to draw these pictures vary; the most na{\"\i}ve work by iterating
the center of a pixel to determine if it lies in the Julia set. 
Milnor's distance estimator algorithm \cite{Mil},
uses classical complex analysis to give a one-pixel estimate of
the Julia set. This algorithm and its modifications work quite well
for many examples, but it is well known that in some particular cases
computation time will grow very rapidly with increase of the resolution.
Moreover, there are examples, even in the family of quadratic polynomials,
when no satisfactory pictures of the Julia set exist.
In this paper we study computability properties of Julia sets of quadratic polynomials.
Under the definition we use, a set is computable, if, roughly speaking,
its image can be generated by a computer with an arbitrary precision.
Under this notion of computability we show:

\medskip
\noindent
{\bf Main Theorem}
{\it There exists a parameter value $c\in\CC$ such that the Julia set of the
quadratic polynomial $f_c(z)=z^2+c$ is not computable.}

\medskip
\noindent
The structure of the paper is as follows. In the Introduction we 
discuss the question of computability of real sets, and make the
relevant definitions. Further in this section we briefly introduce
the reader to the main concepts of Complex Dynamics, and discuss
the properties of Julia sets relevant to us. In the end of the 
Introduction, we outline the conceptual idea of the proof of Main
Theorem. Section \S 3 contains the technical lemmas on which the
argument is based. In \S 4 we complete the proof.

\medskip
\noindent
{\bf Acknowledgement.} 
We are grateful to John Milnor for his encouragement, and for many
helpful suggestions on improving the exposition.
We wish to thank Arnaud Ch{\'e}ritat for
helpful comments on an earlier draft of this paper. 
The first author would like to thank Stephen Cook for
 many discussions on computability of real sets. 

\section{Introduction}
\label{section-intro}

\subsection{Introduction to the Computability of Real Sets}

\subsection*{Classical Computability.} The computability theory in general 
allows us to classify problems into the tractable (``computable")
and intractable (``uncomputable"). All the common computational 
tasks such as integer operations, list sorting, etc. are easily
seen to be computable. On the other hand, there are many uncomputable
problems. 

In the formal setting for the study of computability theory computations 
are performed by objects called the Turing Machines. Turing Machines were 
introduced in 1936 by Alan Turing (see \cite{Tur}) and are accepted by the 
scientific community as the standard model of computation. 
The Turing Machine (TM in short) is capable of solving exactly the same 
problems as an ordinary computer. Most of the time, one can think of 
the TM as a computer program written in any programming language. 
It is important to mention that there are only countably many TMs, which 
can be enumerated in a natural way. 
See \cite{Sip} for a formal discussion on TMs. We define {\em 
computability} as follows.

\begin{defn}
\label{defcomp}
We say that a function $f:\{ 0, 1\}^{*} \rightarrow \{ 0,1\}^{*}$ is 
{\em computable} if there is a TM, which on input string $s$ outputs the 
string $f(s)$. 

We say that the set $L \subset \{ 0,1 \}^{*}$ is {\em computable} or
{\em decidable} if its characteristic function $\chi_L : \{ 0, 1\}^{*} 
\rightarrow \{ 0,1\}$ is computable.
\end{defn}

While most ``common" functions are computable, there are uncountably many 
uncomputable functions and undecidable sets. The best-known intractable 
problems are the {\it Halting Problem} and the {\it solvability of a 
Diophantine equation} (Hilbert's 10-th problem), see \cite{Sip} and 
\cite{Mat} for more information. 

\subsection*{Computability of Real Functions and Sets.} In the present 
paper we are interested in the computability of 
functions $f: \RR^n \rightarrow \RR$ and subsets of 
$\RR^n$, particularly subsets of $\RR^2 \cong \CC$. We cannot directly 
apply Definition \ref{defcomp}   here, since a real number cannot be 
represented in general by finite sequences of bits. 

Denote by $\DD$ the set of the {\it dyadic rationals}, that is, rationals of 
the form $\frac{p}{2^m}$. 
We say that $\phi: \NN \rightarrow \DD$ is an {\it oracle}
for a real number $x$, if $| x - \phi(n)|<2^{-n}$ for all $n \in \NN$. 
In other words, $\phi$ provides a good dyadic approximation for $x$. 
We say that a TM $M^{\phi}$ is an {\it oracle machine}, if at any step of 
the computation $M$ is allowed to query the value $\phi(n)$ for any $n$.
This definition allows us to define the computability of real functions 
on compact sets. 

\begin{defn}
\label{funcomp}
We say that a function $f:[a,b] \rightarrow [c,d]$ is computable, if 
there exits an oracle TM $M^{\phi} (m)$ such that if $\phi$ is an oracle 
for $x \in [a,b]$, then on input $m$, $M^{\phi}$ outputs a $y \in \DD$ 
such that $| y - f(x)|<2^{-m}$.
\end{defn}

This definition was first introduced by Grzegorczyk \cite{Grz} and
Lacombe \cite{Lac}, and follows in the tradition of Computable Analysis
originated by Banach and Mazur in 1937 (see \cite{Maz}).

\noindent
To understand the definition better, the reader without a Computer Science
background should think of a computer program with an instruction
$$\text{ READ real number }x\text{ WITH PRECISION }n(m).$$
On the execution of this command, a dyadic rational $d$ is input from the 
keyboard. This number must not differ from $x$ by more than  $2^{-n(m)}$ 
(but otherwise can be arbitrary). The algorithm then outputs $f(x)$ to precision
$2^{-n}$.

\noindent
In other words, with an access to arbitrarily good approximations for 
$x$, $M$ should be able to produce an arbitrarily good approximation for 
$f(x)$. This definition trivially generalizes to domains of higher 
dimension. See \cite{Ko} for more details. One of the most important 
properties of computable functions is that 

\begin{prop}
\label{cont1}
Computable functions are 
continuous. 
\end{prop}

\noindent
Let $K \subset \RR^k$ be a compact set. We would like to give a definition
for $K$ being computable. Note that saying by a na{\"\i}ve analogy with Definition 
\ref{defcomp}, that $K$ is computable if and only if the 
characteristic function $\chi_K$ is computable does not work here, since by the 
above proposition only continuous functions can be computable.

\noindent
We say that a TM $M$ computes the set $K$ if it approximates $K$ in the {\it 
Hausdorff metric}. Recall that the Hausdorff metric is a metric on 
compact subsets of $\RR^n$ defined by 
\begin{equation}
\label{hausdorff metric}
d_H ( X, Y) =  \inf \{\epsilon > 0 \;|\; X \subset U_{\epsilon} 
(Y)~~\mbox{and}~~  Y \subset U_{\epsilon}(X)\},
\end{equation}

\noindent
where $U_\eps(S)$ is defined as the union of the set of $\eps$-balls with centers in $S$.

\noindent
We introduce
a class $\cC$ of sets which is dense in metric $d_H$ among the compact 
sets and which has a natural correspondence to binary strings. 
Namely $\cC$ is the set of finite unions of dyadic balls:
$$
\cC= \left\{ \bigcup_{i=1}^n \overline{B(d_i, r_i)}~|~~\mbox{where}~~d_i, 
r_i 
\in \DD \right\}.
$$
The following definition is equivalent to the set computability definition
given in \cite{Wei} (see also \cite{WeiPaper}).

\begin{defn}
\label{setcomp}
We say that a compact set $K \subset \RR^k$ is computable, if there exists a TM 
$M(m)$, such that on an input $m\in\NN$, the machine 
$M(m)$ outputs an encoding of $C_m \in 
\cC$ such that $d_H (K, C_m) < 2^{-m}$. 
\end{defn}

To illustrate the robustness of this definition we present the following 
two equivalent characterizations of computable sets. The first one 
relates the definition to computer graphics. It is not stated precisely 
here, but it can be easily made precise. The second one relates the 
computability of sets to the computability of functions as per Definition 
\ref{funcomp}.

\begin{thm}
For a compact $K \subset \RR^k$ the following are equivalent:

(1) $K$ is computable as per Definition \ref{setcomp},

(2) (in the case $k=1, 2$) $K$ can be drawn on a computer screen 
with arbitrarily good precision,

(3) the {\em distance function} $d_K (x) = \inf \{ |x-y|~~|~~y\in K \}$ 
is 
computable as per Definition \ref{funcomp}.
\end{thm}

In the present paper we are interested in   questions concerning 
the computability of the Julia set $J_c = J(f_c) = J(z^2+c)$ (see the next section for the definition). 
Since there 
are uncountably many possible parameter values for $c$, and only countably many TMs, we cannot 
expect for each $c$ to have a machine $M$ such that $M$ computes 
$J_c$. On the other hand, 
it is reasonable to want $M$ to compute $J_c$ with an oracle 
access to $c$. Define the function $J: \CC \rightarrow K^{*}$ ($K^{*}$ is 
the set of all compact subsets of $\CC$) by $J(c)=J(f_c)$. 
In a complete analogy to Definition \ref{funcomp} we can define

\begin{defn}
\label{funcomp2}
We say that a function $f:S \rightarrow K^{*}$ for some bounded set $S\subset \RR^k$ is 
computable, if there exits an oracle TM $M^{\phi} (m)$ such that if 
$\phi$ is an oracle 
for $x \in S$, then on input $m$, $M^{\phi}$ outputs a set $C_m \in \cC$ 
such that $d_H (C_m, f(x))<2^{-m}$.
\end{defn}

\noindent
In the case of Julia sets:

\begin{defn}
\label{Jcomp}
We say that $J_c$ is computable if the function $J: d \mapsto J_d$ is 
computable on the set $\{ c \}$.
\end{defn}

\noindent
The following has been shown (see \cite{thesis}, \cite{Ret}):

\begin{thm}
Denote by $\cH$ the set of parameters $c$ for which $J_c$ is {\em 
hyperbolic}, then 

(i) $J_c$ is computable for all $c \in \cH$, moreover 

(ii) the function $J$ is computable on each bounded subset of $\cH$.
\end{thm}

\noindent
Our goal in this paper is to show that there are values of $c$ for 
which $J_c$ is not computable under Definition \ref{Jcomp}, which is 
the weakest possible definition in this setting. We will be using 
the following version of Theorem \ref{cont1} for set functions.

\begin{thm}
\label{cont2}
Suppose that a TM $M^{\phi}$  computes the function $J$ on a set $S$.
Then $J$ is continuous on $S$ in Hausdorff sense.
\end{thm}
\begin{proof}
Let $c$ be any point in $S$, and let $\ve = 2^{-k}$ be given. Let $\phi$ be 
an oracle for $c$ such that $|\phi(n)-c|<2^{-(n+1)}$ for all $k$. We run 
$M^\phi (k+1)$ with this oracle $\phi$. By the definition of $J$, it outputs a set $L$ which
is a $2^{-(k+1)}$ 
approximation of $J_c$ in the Hausdorff metric. 

The computation is performed in a finite amount of time. Hence 
there is an $m$ such that $\phi$ is only queried with parameters not exceeding $m$. 
Then for any $x$ such that $|x-c|<2^{-(m+1)}$, $\phi$ is a valid oracle for $x$ 
up to parameter value of $m$. In particular, we can create an oracle $\psi$ for $x$ 
that agrees with $\phi$ on $1, 2, \ldots, m$. If $x \in S$, then the execution 
of $M^{\psi}(k+1)$ will be identical to the execution of $M^{\phi}(k+1)$, and it 
will output $L$ which has to be an approximation of $J_x$. Thus we have
$$
d_H(J_c, J_x) \le d_H (J_c, L)+d_H (J_x, L) < 2^{-(k+1)}+2^{-(k+1)} = 2^{-k}.
$$
This is true for any $x \in B(c,2^{-(m+1)})\cap S$. Hence $J$ is continuous
on $S$. 
\end{proof}

\noindent
In the next section we proceed to define Julia sets of rational maps and 
review their basic properties. In particular, towards the end of the introduction,
we will see a mechanism by which the continuity required by \thmref{cont2} may fail.

It should be noted that the question of computability of dynamically generated
fractal sets, such as Julia sets, has been discussed by Blum, Cucker, Shub, and Smale
in \cite{BCSS}. The definition of set computability used in \cite{BCSS} is,
however, quite different from Definition \ref{setcomp}. 
The BCSS model allows infinite-precision arithmetic, but requires
completely accurate pictures to be generated. Under this definition 
all Julia sets but the most trivial ones can be shown to be uncomputable.

\subsection{Julia sets of polynomial mappings}
We recall the main definitions of complex dynamics relevant to our result only briefly;
a good general reference is the book of Milnor \cite{Mil}.
For a rational mapping $R$ of degree $\deg R=d\geq 2$ considered as a dynamical system on the Riemann
sphere
$$R:\riem\to\riem$$
the Julia set is defined as the complement of the set where the dynamics is Lyapunov-stable:

\begin{defn}
Denote $F(R)$ the set of points $z\in\riem$ having an open neighborhood $U(z)$ on which the
family of iterates $R^n|_{U(z)}$ is equicontinuous. The set $F(R)$ is called the Fatou set of $R$
and its complement $J(R)=\riem\setminus F(R)$ is the Julia set.
\end{defn}

\noindent
In the case when the rational mapping is a polynomial $$P(z)=a_0+a_1z+\cdots+a_dz^d:\CC\to\CC$$ an equivalent
way of defining the Julia set is as follows. Obviously, there exists a neighborhood of $\infty$ on $\riem$
on which the iterates of $P$ uniformly converge to $\infty$. Denoting $A(\infty)$ the maximal such domain of attraction
of $\infty$ we have $A(\infty)\subset F(R)$. We then have 
$$J(P)=\partial A(\infty).$$
The bounded set $\riem \setminus \cl A(\infty)$ is called {\it the filled Julia set}, and denoted $K(P)$;
it consists of points whose orbits under $P$ remain bounded:
$$K(P)=\{z\in\riem|\;\sup_n|P^n(z)|<\infty\}.$$

\noindent
For future reference, let us list in a proposition below the main properties of Julia sets:

\begin{prop}
\label{properties-Julia}
Let $R:\riem\to\riem$ be a rational function. Then the following properties hold:
\begin{itemize}
\item $J(R)$ is a non-empty compact subset of $\riem$ which is completely
invariant: $R^{-1}(J(R))=J(R)$;
\item $J(R)=J(R^n)$ for all $n\in\NN$;
\item $J(R)$ is perfect;
\item if $J(R)$ has non-empty interior, then it is the whole of $\riem$;
\item let $U\subset\riem$ be any open set with $U\cap J(R)\neq \emptyset$. Then there exists $n\in\NN$ such that
$R^n(U)\supset J(R)$;
\item periodic orbits of $R$ are dense in $J(R)$.
\end{itemize}
\end{prop}

\noindent
Let us further comment on the last property. For a periodic point $z_0=R^p(z_0)$
of period $p$ its {\it multiplier} is the quantity $\lambda=\lambda(z_0)=DR^p(z_0)$.
We may speak of the multiplier of a periodic cycle, as it is the same for all points
in the cycle by the Chain Rule. In the case when $|\lambda|\neq 1$, the dynamics
in a sufficiently small neighborhood of the cycle is governed by the Mean
Value Theorem: when $|\lambda|<1$, the cycle is {\it attracting} ({\it super-attracting}
if $\lambda=0$), if $|\lambda|>1$ it is {\it repelling}.
Both in the attracting and repelling cases, the dynamics can be locally linearized:
\begin{equation}
\label{linearization-equation}
\psi(R^p(z))=\lambda\cdot\psi(z)
\end{equation}
where $\psi$ is a conformal mapping of a small neighborhood of $z_0$ to a disk around $0$.
By a classical result of Fatou, a rational mapping has at most finitely many non-repelling
periodic orbits. Therefore, we may refine the last statement of \propref{properties-Julia}:

\begin{itemize}
\item {\it repelling periodic orbits are dense in $J(R)$}.
\end{itemize}

\noindent
In the case when $|\lambda|=1$, so that $\lambda=e^{2\pi i\theta}$, $\theta\in\RR$, 
 the simplest to study is the {\it parabolic case} when $\theta=n/m\in\QQ$, so $\lambda$ 
is a root of unity. In this case $R^p$ is not locally linearizable; it is not hard to see that $z_0\in J(R)$.
 In the complementary situation, two non-vacuous possibilities  are considered:
{\it Cremer case}, when $R^p$ is not linearizable, and {\it Siegel case}, when it is.
In the latter case, the linearizing map $\psi$ from (\ref{linearization-equation}) conjugates
the dynamics of $R^p$ on a neighborhood $U(z_0)$ to the irrational rotation by angle $\theta$
(the {\it rotation angle})
on a disk around the origin. The maximal such neighborhood of $z_0$ is called a {\it Siegel disk}.
Siegel disks will prove crucial to our study, and will be discussed in more detail in the 
next section. 

To conclude the discussion of the basic properties of Julia sets, let us consider the simplest
examples of non-linear rational endomorphisms of the Riemann sphere, the quadratic polynomials.
Every affine conjugacy class of quadratic polynomials has a unique representative of the
form $f_c(z)=z^2+c$, the family
$$f_c(z)=z^2+c,\;c\in\CC$$
is often referred to as {\it the quadratic family}.
For a quadratic map the structure of the Julia set is governed by the behavior of the orbit of the only
finite critical point $0$. In particular, the following dichotomy holds:

\begin{prop}
\label{quadratic-Julia}
Let $K=K(f_c)$ denote the filled Julia set of $f_c$, and $J=J(f_c)=\partial K$. Then:
\begin{itemize}
\item $0\in K$ implies that $K$ is a connected, compact subset of the plane with connected complement;
\item $0\notin K$ implies that $K=J$ is a planar Cantor set.
\end{itemize}
\end{prop}

\noindent
The {\it Mandelbrot set} $\cM\subset \CC$ is defined as the set of parameter values $c$ for which 
$J(f_c)$ is connected.

\subsection*{Continuity of the dependence $c\mapsto J(f_c)$}
A natural question to pose for polynomials in the quadratic family is whether the Julia set varies
continuously with the parameter $c$. To make sense of this question, recall the definition of the
Hausdorff distance $\dist_H$ between compact sets $X$, $Y$ in the plane (\ref{hausdorff metric}).
It turns out that the dependence $c\mapsto J(f_c)$ is discontinuous in the Hausdorff distance.
For an excellent survey of this problem see the paper of Douady \cite{Do}. The discontinuity which
has found most interesting dynamical applications occurs at parameter values for which $f_c$ has
a parabolic point. We, however, will employ a more obvious discontinuity which is related to Siegel
disks. Let us first note that by a result of Douady and Hubbard \cite{orsay-notes} a quadratic polynomial 
has at most one non-repelling cycle
in $\CC$. In particular, there is at most one cycle of Siegel disks.

\begin{prop} Let $c_*\in\cM$ be a parameter value for which $f_c$ has a Siegel disk. Then the map
$c\mapsto J(f_c)$ is discontinuous at $c_*$.
\end{prop}

\begin{pf}
Let $z_0$ be a Siegel periodic point of $f_c$ and denote $\Delta$ the Siegel disk around $\zeta_0$, 
$p$ its period, and $\theta$ the rotation angle.
By the Implicit Function Theorem, there exists a holomorphic mapping $\zeta:U(c_*)\to \CC$ such that
$\zeta(c_*)=z_0$ and $\zeta(c)$ is fixed under $(f_c)^p$. The mapping 
$$\nu:c\mapsto D(f_c)^p(\zeta(c))$$
is holomorphic, hence it is either constant or open. If it is constant, all quadratic polynomials have a Siegel
disk. This is not possible: for instance, $f_{1/4}$ has a parabolic fixed point, and thus no other
non-repelling cycles. Therefore, $\nu$ is open, and in particular, there is a sequence of 
parameters $c_n\to c_*$ such that $\zeta(c_n)$ has multiplier $e^{2\pi i p_n/q_n}$. 
Since $\zeta(c_n)$ is parabolic, it lies in the Julia set of $f_{c_n}$. Hence
$$\dist_H(J(f_{c_*}),J(f_{c_n}))>\dist(c_*,\partial\Delta)/2$$
for $n$ large enough.
\end{pf}

\noindent
Thus an arbitrarily small change of the multiplier of the Siegel point may lead to an implosion of the 
Siegel disk -- its inner radius collapses to zero.
We make a note of an immediate consequence of the above proposition and
\thmref{cont2}:

\begin{prop}
For any TM $M^\phi(n)$ with an oracle for $c\in\CC$ denote $S_M$ the set
of all values of $c$ for which $M^\phi$ computes $J_c$. Then $S_M\neq \CC$.
\end{prop}

\noindent
In other words, a single algorithm for computing all quadratic Julia sets
does not exist.

\subsection*{Siegel disks of quadratic maps}
Let us discuss in more detail the occurrence of Siegel disks in the quadratic family.
For a number $\theta\in [0,1)$ denote $[r_1,r_2,\ldots,r_n,\ldots]$, $r_i\in\NN\cup\{\infty\}$ its possibly finite 
continued fraction expansion:
\begin{equation}
\label{cfrac}
[r_1,r_2,\ldots,r_n,\ldots]\equiv\cfrac{1}{r_1+\cfrac{1}{r_2+\cfrac{1}{\cdots+\cfrac{1}{r_n+\cdots}}}}
\end{equation}
Such an expansion is defined uniquely if and only if $\theta\notin\QQ$. In this case, the {\it rational 
convergents } $p_n/q_n=[r_1,\ldots,r_{n}]$ are the closest rational approximants of $\theta$ among the
numbers with denominators not exceeding $q_n$. In fact, setting $\lambda=e^{2\pi i\theta}$, we have
$$|\lambda^h-1|>|\lambda^{q_n}-1|\text{ for all }0<h<q_{n+1},\; h\neq q_n.$$
The difference $|\lambda^{q_n}-1|$ lies between $2/q_{n+1}$ and $2\pi/q_{n+1}$,
therefore the rate of growth of the denominators $q_n$ describes how well 
$\theta$ may be approximated with rationals.

\begin{defn}
The {\it diophantine numbers
of order k}, denoted $\cD(k)$
is the following class of irrationals ``badly'' approximated by rationals.
 By definition, $\theta\in\cD(k)$ if there exists $c>0$ such that
$$q_{n+1}<cq_n^{k-1}$$
\end{defn}

\noindent
The numbers $q_n$ can
be calculated from the recurrent relation
$$q_{n+1}=r_{n+1}q_n+q_{n-1},\text{ with }q_0=0,\; q_1=1.$$ Therefore, $\theta\in\cD(2)$ if and only if the sequence $\{r_i\}$
is bounded. Dynamicists call such numbers {\it bounded type} (number-theorists prefer {\it constant type}). 
An extreme example of a number of bounded type is the golden mean
$$\theta_*=\frac{\sqrt{5}-1}{2}=[1,1,1,\ldots].$$
The set $$\displaystyle\cD(2+)\equiv\bigcap_{k>2}\cD_k$$
has full measure in the interval $[0,1)$. In 1942 Siegel showed:

\begin{thm}[\cite{siegel}]
Let $R$ be an analytic map with an periodic point $z_0\in\riem$ of period $p$. Suppose the multiplier of the cycle
$$\lambda=e^{2\pi i\theta}\text{ with }\theta\in\cD(2+),$$
then the local linearization equation (\ref{linearization-equation}) holds.
\end{thm}

\noindent
The strongest known generalization of this result was proved by Brjuno in 1972:
\begin{thm}[\cite{Bru}]
Suppose
\begin{equation}
\label{brjuno}
B(\theta)=\displaystyle\sum_n\frac{\log(q_{n+1})}{q_n}<\infty,
\end{equation}
then the conclusion of Siegel's Theorem holds.
\end{thm}

\noindent
Note that a quadratic polynomial with a fixed Sigel disk with rotation angle $\theta$ after an affine
change of coordinates can be written as 
\begin{equation}
\label{form-1}
P_\theta(z)=z^2+e^{2\pi i \theta}z.
\end{equation}
\noindent
In 1987 Yoccoz \cite{Yoc} proved the following converse to Brjuno's Theorem:

\begin{thm}[\cite{Yoc}]
Suppose that for $\theta\in[0,1)$ the polynomial $P_\theta$ has a Siegel point at the origin.
Then $B(\theta)<\infty$.
\end{thm}

\noindent
The numbers satisfying (\ref{brjuno}) are called Brjuno numbers; the set of all Brjuno numbers will be denoted $\cB$.
It is evident that $\cup\cD(k)\subset \cB$. The sum of the series (\ref{brjuno}) is called the Brjuno function. 
For us a different characterization of $\cB$ will be more useful. Inductively define $\theta_1=\theta$
and $\theta_{n+1}=\{1/\theta_n\}$. In this way, 
$$\theta_n=[r_{n},r_{n+1},r_{n+2},\ldots].$$
We define the {\it Yoccoz's Brjuno function} as
$$\Phi(\theta)=\displaystyle\sum_{n=1}^{\infty}\theta_1\theta_2\cdots\theta_{n-1}\log\frac{1}{\theta_n}.$$
One can verify that $$B(\theta)<\infty\Leftrightarrow \Phi(\theta)<\infty.$$
The value of the function $\Phi$ is related to the size of the Siegel disk in the following way.

\begin{defn}
Let $P(\theta)$ be a quadratic polynomial with a Siegel disk $\Delta_\theta\ni 0$. Consider a conformal isomorphism
$\phi:\DD\mapsto\Delta$ fixing $0$. The {\it conformal radius of the Siegel disk $\Delta_\theta$} is
the quantity
$$r(\theta)=|\phi'(0)|.$$
For all other $\theta\in[0,\infty)$ we set $r(\theta)=0$. 
\end{defn} 

\noindent
By the Koebe One-Quarter Theorem of classical complex analysis, the internal radius of $\Delta_\theta$ is at least
$r(\theta)/4$. Yoccoz \cite{Yoc} has shown that the sum 
$$\Phi(\theta)+\log r(\theta)$$
is bounded from below independently of $\theta\in\cB$. Recently, Buff and Ch{\'e}ritat have greatly improved this result
by showing that:

\begin{thm}[\cite{BC}]
\label{phi-cont}
The function $\theta\mapsto \Phi(\theta)+\log r(\theta)$ extends to $\RR$ as a 1-periodic continuous
function.
\end{thm}

\noindent
We remark that the following stronger conjecture exists (see \cite{MMY}):

\medskip
\noindent
{\bf Marmi-Moussa-Yoccoz Conjecture.} \cite{MMY} {\it The function $\theta\mapsto \Phi(\theta)+\log r(\theta)$ is H{\"o}lder of exponent $1/2$.}

\subsection*{Dependence of the conformal radius of a Siegel disk on the parameter}
In this section we will show that the conformal radius of a Siegel disk varies continuously with the Julia set.
To that end we will need a preliminary definition:

\begin{defn}
Let $(U_n,u_n)$ be a sequence of topological disks $U_n\subset\CC$ with marked points $u_n\in U_n$.
The {\it kernel} or {\it Carath{\'e}odory} convergence $(U_n,u_n)\to (U,u)$ means the following:
\begin{itemize}
\item $u_n\to u$;
\item for any compact $K\subset U$ and for all $n$ sufficiently large, $K\subset U_n$;
\item for any open connected set $W\ni u$, if $W\subset U_n$ for infinitely many $n$, then $W\subset U$.
\end{itemize}
\end{defn}

\noindent
The topology on the set of pointed domains which corresponds to the above definition of convergence is again
called {\it kernel} or {\it Carath{\'e}odory} topology. The meaning of this topology is as follows.
For a pointed domain $(U,u)$ denote 
$$\phi_{(U,u)}:\DD\to U$$
the unique conformal isomorphism with $\phi_{(U,u)}(0)=u$, and $(\phi_{(U,u)})'(0)>0$.  
We again denote $r(U,u)=|(\phi_{(U,u)})'(0)|$ the conformal radius of $U$ with respect to $u$.

By the Riemann Mapping
Theorem, the correspondence $$\iota:(U,u)\mapsto \phi_{(U,u)}$$
establishes a bijection between marked topological disks properly contained in $\CC$ and univalent maps $\phi:\DD\to\CC$
with $\phi'(0)>0$.
The following theorem is due to Carath{\'e}odory, a proof may be found in  \cite{Pom}:

\begin{thm}[{\bf Carath{\'e}odory Kernel Theorem}]
The mapping $\iota$ is a homeomorphism with respect to the Carath{\'e}odory topology on domains and the
compact-open topology on maps.
\end{thm}

\noindent
\begin{prop}
\label{radius-continuous}
The conformal radius of a quadratic Siegel disk varies continuously with respect to the Hausdorff 
distance on Julia sets.
\end{prop}

\begin{pf}
To fix the ideas, consider the family $P_\theta$ with $\theta\in\cB$ and denote $\Delta_\theta$ the Siegel disk
of $P_\theta$. It is easy to see that the Hausdorff convergence $J(P_{\theta_n})\to J(P_\theta)$ implies the
Carath{\'e}odory convergence of the pointed domains
$$(\Delta_{\theta_n},0)\to(\Delta,0).$$
The proposition follows from this and the Carath{\'e}odory Kernel Theorem.
\end{pf}

\noindent
In fact, we can state the following quantitative version of the above result. For the proof, based on
Koebe Theorem, see e.g. \cite{RZ}:
\begin{lem}
\label{variation conf radius}
Let $U$ be a simply-connected  bounded subdomain of $\CC$ containing the point $0$ in the interior.
Suppose $V\subset U$ is a simply-connected subdomain of $U$, and $\partial V\subset U_\eps(\partial U)$.
Then 
$$0<r(U,0)-r(V,0)\leq 4\sqrt{r(U,0)}\sqrt{\eps}.$$

\end{lem}

\noindent
An immediate corollary is:
\begin{cor}
\label{noncomp radius}
Suppose the function $r(\theta)$ is uncomputable on the set $\{\theta_0\}$. Then the function
$\theta\mapsto J(P_\theta)$ is also uncomputable at the same point.

\end{cor}
\begin{pf}
Assume that $J(P_{\theta_0})$ is computable. Using the output of the TM computing this Julia set
in an obvious way,
for each $\eps>0$ we can obtain a  domain $V\in\cC$ such that 
$$V\subset \Delta_{\theta_0}\text{ and }
d_H(\partial V,\partial \Delta_{\theta_0})<\eps.$$ 
By Schwarz Lemma, the conformal radius $r(\theta_0)<2$. Hence, 
by \lemref{variation conf radius}, 
$$|r(V,0)-r(\theta_0)|<\delta=8\sqrt{\eps}.$$
Using any constructive version of the Riemann Mapping Theorem (see e.g. \cite{BB}), we 
can compute $r(V,0)$ to precision $\delta$, and hence know $r(\theta_0)$ up to an error of $2\delta$.
Given that $\delta$ can be made arbitrarily small, we have shown that $r(\theta_0)$ is computable.

\end{pf}

\noindent
We also state for future reference the following proposition:

\begin{prop}
\label{r doest drop}
Let $\{\theta_i\}$ be a sequence of Brjuno numbers such that $\theta_i\to\theta$ and
$\overline{\lim}\; r(\theta_i)=l>0$. Then $\theta$ is also a Brjuno number and $r(\theta)\geq l$.
\end{prop}
\begin{pf}
Denote $\phi_i\equiv \phi_{(\Delta_{\theta_i},0)}.$ 
Note that by Schwarz Lemma, the inverse $\psi_i\equiv (\phi_i)^{-1}$ linearizes $P_{\theta_i}$
on $\Delta_{\theta_i}$. By passing to a subsequence we can assure that 
$\phi_i\to\phi$ locally uniformly, and $\phi'(0)\geq l$. By continuity, $\phi^{-1}$ is a linearizing
coordinate for $P_\theta$, so $\theta$ is a Brjuno number. Moreover, $\phi(\DD)\subset \Delta_\theta$,
and so by Schwarz Lemma $r(\theta)\geq l$.
\end{pf}

\subsection*{Non-computability of the Yoccoz's Brjuno function.}
In addition to the non-compu{\-}ta{\-}bi{\-}lity of the conformal radius, 
we also prove a  non-computability 
result for the Yoccoz's 
Brjuno function $\Phi$:

\begin{thm}[{\bf Non-computability of $\mathbf \Phi$}]
\label{secondmain}
There exists a  parameter value $\theta\in\RR/\ZZ$ such that 
$\Phi(\theta)<\infty$, and $\Phi(\theta)$ is not computable by any Turing Machine with an oracle for $\theta$.
\end{thm}

\noindent
It is worth noting that Marmi-Moussa-Yoccoz Conjecture as stated above and 
 Theorem \ref{secondmain} imply the existence of a non-computable quadratic Julia set.
To see this, we first formulate:

\begin{tet}
If Marmi-Moussa-Yoccoz Conjecture holds, then the function 
$$\upsilon: \theta\mapsto \Phi(\theta)+\log r(\theta)$$
is computable by one Turing Machine on the entire interval $[0,1]$. 
\end{tet}

\noindent
We use the following result of Buff and Ch{\'e}ritat (\cite{BC}). 

\begin{lem}[\cite{BC}]
For any rational point $\theta = \frac{p}{q} \in [0,1]$ denote, as before,
$$P_{\theta}(z)= 
e^{2 \pi i \theta} z + z^2,$$ and let the Taylor expansion of $P_\theta^{\circ q} (z)$ at $0$ start with  
$$
P_\theta^{\circ q} (z) = z + A z^{q+1} + \ldots,\text{ for }q\in\NN
$$
Let $L(\theta) = \left( \frac{1}{q A}\right)^{1/q}$. Denote by $\Phi_{trunc}$
the modification of $\Phi$ applied to rational numbers where the sum is truncated
before the infinite term. Then we have the following explicit formula
for computing $\upsilon(\theta)$:
\beq
\label{expl-rat}
\upsilon(\theta) = \Phi_{trunc} (\theta) + \log L(\theta) + \frac{\log 2\pi}{q}.
\eeq
\end{lem}

\noindent
Equation \eref{expl-rat} allows us to compute the value of $\upsilon$ easily at
every rational $\theta \in \QQ \cap [0,1]$ with an arbitrarily good precision. 
In addition, assuming the conjecture, we have $| \upsilon(x) - \upsilon(y)|< 2^{-n}$
whenever $|x-y| < c \cdot 2^{-2 n}$ for some constant $c$, hence $\upsilon$ has 
an (easily) computable modulus of continuity. These two facts together
imply that $\upsilon$ is computable by a single machine of the interval $[0,1]$ 
(see for example Proposition 2.6 in \cite{Kosurv}). This implies the Conditional
Implication.

\noindent
The following conditional result follows:

\begin{lem}
[{\bf Conditional}] \label{2main-equiv}
Suppose the Conditional Implication holds. Let $\theta \in [0,1]$
be such that $\Phi(\theta)$ is finite. Then there is an oracle Turing 
Machine $M^{\phi}_1$ computing $\Phi(\theta)$ with an oracle access to 
$\theta$ if and only if there is an oracle Turing 
Machine $M^{\phi}_2$ computing $r(\theta)$ with an oracle access to 
$\theta$. 
\end{lem}

\begin{proof}
Suppose that $M^{\phi}_1$ computes $\Phi(\theta)$ for some 
$\theta$. Let $M^{\phi}$ be the machine uniformly computing the 
function $\upsilon$. Then we can use $M_1^{\phi}$ and $M^{\phi}$ 
to compute $\log r(\theta) = \upsilon(\theta) - \Phi (\theta)$ with 
an arbitrarily good precision. We can then use this construction 
to give a machine $M_2^{\phi}$ which computes $r(\theta)$. 

The opposite direction is proved analogously. 
\end{proof}

\noindent
Lemma \ref{2main-equiv} with Theorem \ref{secondmain} imply that 
there is a $\theta$ for which $r(\theta)$ is non-computable.
\corref{noncomp radius} implies
that for this value of $\theta$ the Julia
set of $P_\theta$ is also non-computable.

\noindent
Note that for the proof of Conditional Implication we did not need
the full power of the conjecture. All we needed is some {\em computable}
bound on the modulus of continuity of $\upsilon$. 

\subsection*{Outline of the construction of a non-computable quadratic Julia set}
We are now prepared to outline the idea of our construction. 
The outline given below is rather rough and suffers from obvious logical 
deficiencies. However, it captures the idea of the proof in a simple to understand form.
Suppose that every Julia set of a polynomial $P_\theta$ is computable by 
an oracle machine $M^{\phi}$, where $\phi$ represents $\theta$. There are 
countably many machines, so we can enumerate them $M^{\phi}_1, M^{\phi}_2, 
\ldots$. Denote by $S_i$ the domain on which $M^{\phi}_i$ computes 
$J_{P_\theta}$ 
properly. Then we must have:
\begin{enumerate}
\item 
$\CC = \bigcup_{i=1}^\infty S_i,$

\item 
for each $i$, the function $J: \theta \mapsto J({P_\theta})$ is continuous on 
$S_i$.
\end{enumerate}

\realfig{aa}{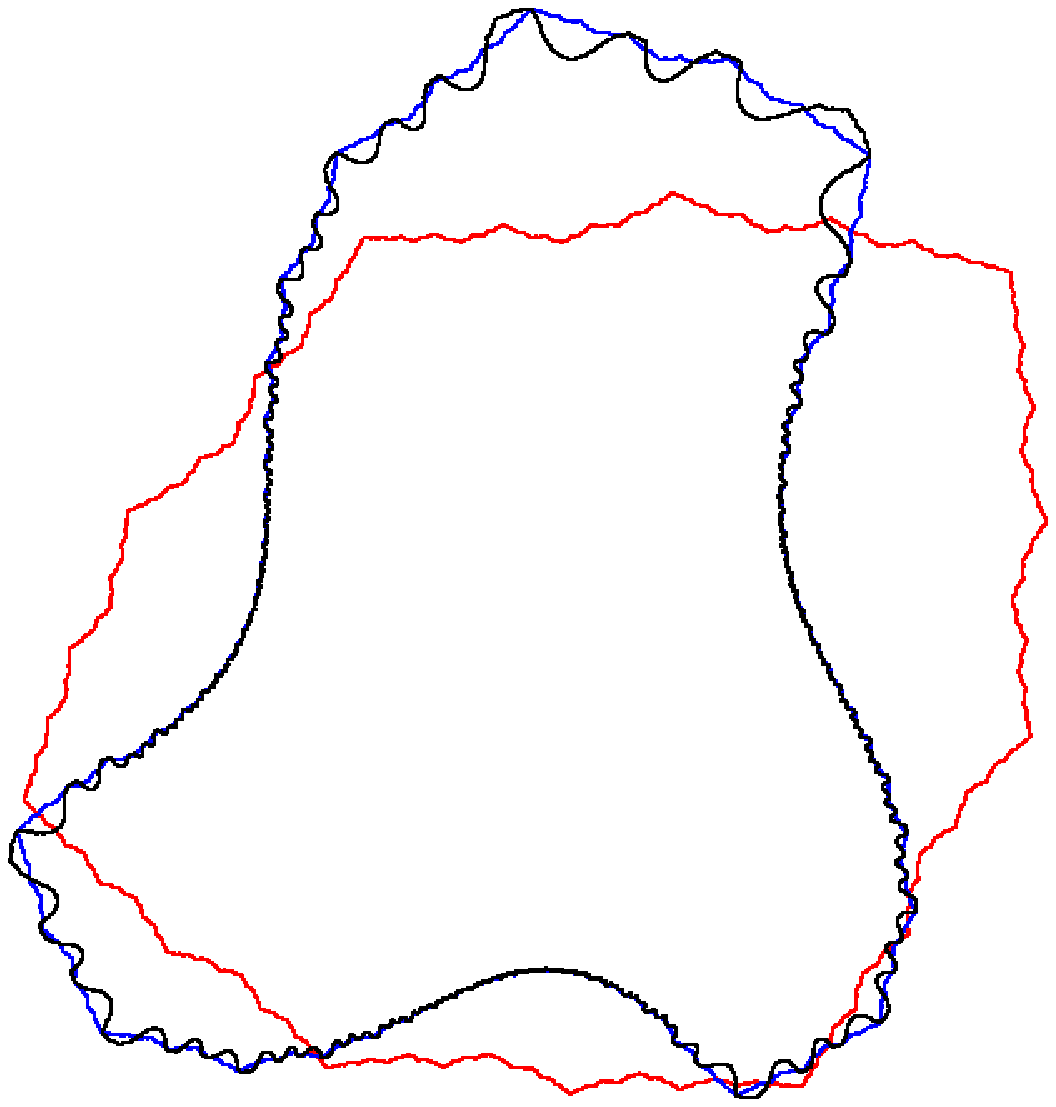}{
The Siegel disks of $P_\theta$ for $\theta$ given by the continued fractions
$[1,1,1,\ldots]$, $[1,1,1,20,1,\ldots]$, and $[1,1,1,20,1,1,1,30,1,\ldots]$}{7cm}

\noindent
Let us start with a machine $M^\phi_{n_1}$ which computes $J(P_{\theta_*})$ for $\theta_*=[1,1,1,\ldots]$.
If any of the digits $r_i$ in this infinite continued fraction is changed to a sufficiently large $N\in\NN$,
the conformal radius of the Siegel disk will become small. For $N\to\infty$ the Siegel disk will implode
and its center will become a parabolic fixed  point in the Julia set.
Given the continuity of the dependence of the conformal radius of the Siegel disk on the Julia set,
we have the following:

\medskip
\noindent 
There exists $i_1>1$ such that for every $\theta_1$ whose continued fraction starts with $i_1$ 
ones, for the Julia set of $P_{\theta_1}$ to be  computable by $M^\phi_{n_1}$, it must possess a Siegel 
disk of a conformal radius 
$r(\theta_1)>r(\theta_*)(1-1/8).$

\medskip
\noindent
We can thus ``fool'' the machine $M_{n_1}^\phi$ by selecting $\theta_1$ given by a continued fraction
where all digits are ones except $r_{i_1}=N_1>>1$. If we are careful, we can do it so that
\beq \label{exam1} r(\theta_*)(1-1/4)<r(\theta_1)<r(\theta_*)(1-1/8).\eeq
To ``fool'' the machine $M_{n_2}^\phi$ we then change a digit $r_{i_2}$ for $i_2>i_1$ sufficiently far in the 
continued fraction of $\theta_1$ to a large $N_2$. In this way, we will obtain a Brjuno number $\theta_2$
for which
\beq \label{exam2} r(\theta_*)(1-1/4-1/8)<r(\theta_2)<r(\theta_*)(1-1/4).\eeq
Continuing in this manner we will arrive at a limiting Brjuno number $\theta_\infty$ for which the Julia
set is uncomputable. 
To make such a scheme work, we need a careful analysis of the dependence of the conformal radius on
the parameter. In this a key role is played by \thmref{phi-cont} of Buff and Ch{\'e}ritat 
which allows us to obtain a controlled change in the value of $r(\alpha)$ by changing $\Phi(\alpha)$.
The relevant analysis is carried out in the next section.

\subsection*{Main analytic result.} We formalize the strategy outlined above as follows:

\begin{thm}
\label{no partition}
There does not exist a partition of the circle $\RR/\ZZ$ into a countable union of 
sets $S_i$ such that for every $i$ the function
$r(\theta)$ restricted to $S_i$ is continuous.
\end{thm}

\noindent
The above formulation was suggested to us by John Milnor. Let us show how the Main Theorem and \thmref{secondmain}
follow from \thmref{no partition}.

\begin{proof}[Proof of Main Theorem, assuming \thmref{no partition}]
First we observe that there exists a parameter $\theta_0\in\RR/\ZZ$ such that the function $r(\theta)$ is
uncomputable on $\{\theta_0\}$. Indeed, assume the contrary. There are only countably many Turing Machines
with an oracle for $\theta$. We enumerate them $M_i$, $i\in\NN$ in some arbitrary way (for instance, using 
the lexicographic order). Denote
$$S_i=\{\theta\in\RR/\ZZ\text{ such that TM }M_i^\phi\text{ computes }r(\theta)\}$$
By \propref{cont1} the function $r(\theta)$ is continuous on each of the $S_i$'s, and we 
arrive at a contradiction with \thmref{no partition}.

Now let us prove Main Theorem, again arguing by contradiction. Assume that for every $c\in\CC$ there
exists a TM $M^\phi$ with an oracle for $c$ which computes $J_c$. Let $P_\theta=z^2+e^{2\pi i\theta}z$
as before. The affine change of coordinates transforming it into an element of the family $f_c$ is 
computable explicitly, and we have $$c=c(\theta)=\lambda^2/4-\lambda/2\text{ where }\lambda=e^{2\pi i\theta}.$$
This implies that we can simulate an oracle for $c$ given an oracle for $\theta$.

Set $c_0=c(\theta_0)$ and consider the oracle TM $M^\phi$ computing the Julia set of $f_{c_0}$.
By the above considerations, there exists an oracle TM $\widetilde{M}^\psi$ with an oracle for $\theta\in\RR/\ZZ$
which computes $J(P_{\theta_0})$. This contradicts \corref{noncomp radius} and the proof is complete.

\end{proof}

\begin{proof}[Proof of \thmref{secondmain} assuming \thmref{no partition}]
Assume the contrary. Again, order in a sequence $M_i^\phi$, $i\in\NN$ all TMs with
an oracle for $\theta\in\RR/\ZZ$. Let 
$$\Omega_i=\{\theta\in\RR/\ZZ\text{ such that }M_i^\phi\text{ computes the value of }\Phi(\theta)\}.$$
Denote $\Omega_0$ the set of all $\theta$ with $\Phi(\theta)=\infty$. The value of $r(\theta)$ on $\Omega_0$
is thus identically $0$.

Denote $$\upsilon(\theta)=\Phi(\theta)+\log r(\theta),$$
which by \cite{BC} continuously extends to $\RR/\ZZ$.
Given \propref{cont1}, the function $\Phi(\theta)$ is continuous on each $\Omega_i$, $i\in\NN$,
and hence so is $$r(\theta)=\exp(\upsilon(\theta)-\Phi(\theta)).$$
By our assumption, $$\RR/\ZZ=\cup_{i=0}^\infty\Omega_i,$$
and we arrive at a contradiction with \thmref{no partition}.

\end{proof}

\subsection*{A note on the connection with \cite{BC2}}
A.~Ch{\'e}ritat has pointed to us that methods of \cite{BC2}, where Siegel disks 
with smooth boundaries are constructed for the quadratic family can be used to
derive the Main Theorem. We discuss this in the section following the proofs
of the main theorems. We note here that the argument we give is based on quite
elementary estimates of the function $\Phi$ and is thus accessible to non-dynamicists.
It has an added advantage of yielding \thmref{secondmain}.

\section{Making Small Changes to $\Phi$}
\label{sec changes}

\subsection{Small Changes to $\Phi$}

A key step of the construction outlined above is 
making careful  adjustments of $r(\theta_i)$ as 
in the first two steps \eref{exam1} and \eref{exam2} above. 
We do not have a direct control over the value of 
$r(\al)$, but Buff and Ch{\'e}ritat's Theorem \ref{phi-cont}
shows that small decreases of $r(\al)$ we would like to 
make correspond to a small controlled increment of the
value of $\Phi(\al)$. Estimates of a similar nature
has appeared in the works of various authors (compare, for example, with \cite{BC2}).

For a number $\gamma=[a_1,a_2,\ldots]\in\RR\setminus\QQ$ we denote
$$
\displaystyle\al_i(\gamma) = \frac{1}{a_i + \displaystyle\frac{1}{a_{i+1} + \displaystyle\frac{1}{a_{i+2}+\ldots}}}, 
$$
so that
$$
\Phi(\gamma) = \sum_{n\ge 1} \al_1(\gamma) \al_2(\gamma) \ldots 
\al_{n-1}(\gamma) \log \frac{1}{\al_n(\gamma)}.
$$
The main goal of this section is 
to prove the following two lemmas:

\begin{lem}
\label{smlchg}
For any initial segment $I = [a_1, a_2, \ldots, a_n]$, write 
$\om = [a_1, a_2, \ldots, a_n, 1, 1, 1, \dots]$. Then 
for any $\ve>0$, there is an $m>0$ and an integer $N$
such that if we write $\be^N = [a_1, a_2, \ldots, a_n, 1, 1, \ldots, 
1, N, 1, 1, \ldots]$, where the $N$ is located in the $n+m$-th 
position, and
$$
\Phi(\om) + \ve < \Phi(\be^N) < \Phi(\om) + 2 \ve. 
$$
\end{lem}

\begin{lem}
\label{notdeclem}
For  $\om$ as above, for any $\ve>0$ there is an $m_0>0$,  
such that for any $m \ge m_0$, and for any tail 
$T = [a_{n+m}, a_{n+m+1}, \ldots]$ if we denote 
$$ \be^{T} = [a_1, a_2, \ldots, a_n, 1, 1, \ldots, 1, a_{n+m}, a_{n+m+1}, 
\ldots],$$ then
$$
\Phi(\be^{T}) > \Phi (\om) - \ve. 
$$
\end{lem}

\noindent
The proof is technical and will require some preparation. 
For Lemma \ref{smlchg}, the idea is to choose an $m$ large enough, so that changing 
$a_{n+m}$ (which will eventually be $N$) by $1$ changes 
 the value of $\Phi$ by a very small amount ($< \ve$). 
When $N \rightarrow \infty$, $\Phi(\om) \rightarrow \infty$, 
hence the value of $\Phi$ must hit the interval $(\Phi(\om)+\ve, 
\Phi(\om)+ 2 \ve)$. 

Denote 
$$
\Phi^{-} (\om) = \Phi(\om) - \al_1 \al_2 \ldots \al_{n+m-1} \log 
\frac{1}{\al_{m+n}}. 
$$
The value of the integer $m>0$ is yet to be determined.
Denote  $$\be^N = [a_1, a_2, \ldots, a_n, 1, 1, \ldots,
1, N, 1, 1, \ldots].$$ We prove the following 

\begin{lem}
\label{lem7} 
For any $N$ and $i\le n+m$ we have $$\left|
\log{\displaystyle\frac{\al_{i}(\be^{N})}{\al_{i}(\be^{N+1})}} \right|
< \displaystyle\frac{2^{i-(n+m)}}{N}.$$
\end{lem}

\begin{proof}

We prove the lemma by induction on $i$, starting from the base case 
$i=n+m$, and proceeding down to $i=0$. The base case is 
$i=n+m$, we want to prove
$$
\left|
\log{\displaystyle\frac{\al_{n+m}(\be^{N})}{\al_{n+m}(\be^{N+1})}} \right|
< \displaystyle\frac{1}{N}.$$
We have, 
$$r_{n+m}=\displaystyle\frac{\al_{n+m}(\be^N)}
{\al_{n+m}(\be^{N+1})} = \displaystyle\frac{\displaystyle\frac{1}{N+ 1/\phi}}{\displaystyle\frac{1}{N+1+1/\phi}} =
\displaystyle\frac{N+1+1/\phi}{N+1/\phi} = 1 + \displaystyle\frac{1}{N+1/\phi},
$$
where $\phi = \left(\sqrt{5}+1\right)/2$. 
Hence
$$
1 <r_{n+m} < 1 + \displaystyle\frac{1}{N} < e^{1/N},
$$
and 
$
\left|\log{r_{n+m}}\right|
< \displaystyle\frac{1}{N}.$

\medskip
\noindent
{\bf Induction step.} Supposing that the statement is true for 
$i+1$, we prove it for $i$. We have 
$$
\frac{\al_{i}(\be^{N})}{\al_{i}(\be^{N+1})} = 
\frac{\displaystyle\frac{1}{a_i+\al_{i+1}(\be^{N})}}{\displaystyle\frac{1}{a_i+\al_{i+1}(\be^{N+1})}}=
\frac{a_i+\al_{i+1}(\be^{N+1})}{a_i+\al_{i+1}(\be^{N})}.
$$
Suppose that $\al_{i+1}(\be^{N+1}) \ge \al_{i+1} (\be^N)$. Then we know that 
$$\frac{\al_{i+1}(\be^{N+1})}{\al_{i+1}(\be^{N})} < e^{2^{i+1-(n+m)}/N},$$ and 
we want to prove that $$\frac{a_i+\al_{i+1}(\be^{N+1})}{a_i+\al_{i+1}(\be^{N})}< 
e^{2^{i-(n+m)}/N},$$ since this expression is obviously bigger than $1$. 
The situation is very similar in the case when $\al_{i+1}(\be^{N+1}) \le 
\al_{i+1} 
(\be^N)$, with the numerator and the denominator exchanged. 

In other words, it is enough to prove that for $0< d < c < 1$ and a pair of integers $r \ge 
0$, $k \ge 1$ and $\al>0$,  $$\frac{c}{d} < e^{\al}\text{ implies that }\frac{k+c}{k+d} <
e^{\al/2}.$$ First of all, it is easy to see that for $k \ge 1$, 
$$\frac{k+c}{k+d} \le \frac{1+c}{1+d},$$ hence it suffices to show that 
$$\frac{1+c}{1+d}<e^{\al/2}=\left( e^{\al}\right)^{1/2}.$$ Thus we need to demonstrate
 that $$\frac{1+c}{1+d} < \left( \frac{c}{d} \right)^{1/2}.$$
This is equivalent to $$(1 + 2 c + c^2) d < ( 1 + 2 d + d^2) c\;\Leftrightarrow\;
 d + c^2 d < c + d^2 c\;\Leftrightarrow\; c d (c - d) < c - d.$$ The last inequality 
holds, since $ c d < 1$ and $c-d >0$. 
\end{proof}

\medskip
\noindent
The following lemma is proven by induction exactly as the previous one with 
a different base. 

\begin{lem}
\label{lem2}
Let $\ga_1$ and $\ga_2$ be two numbers whose continued fraction 
expansions coincide in the first $n+m-1$ terms $[a_1, a_2, \ldots, a_{n+m-1}]$. 
Then for any $i <  n+m$ we have $$\left| 
\log{\displaystyle\frac{\al_{i}(\ga_1)}{\al_{i}(\ga_2)}} \right|
< 2^{i-(n+m)+1}.$$
In particular, this applies with $\ga_1 = \be^N$ and $\ga_2 = \be^1$. 
\end{lem}

\begin{proof}
The proof goes by induction exactly as in Lemma \ref{lem7}. We need to 
verify the base case $i=n+m-1$. For this value of $i$, 
$$\al_{n+m-1}(\ga_1) = \frac{1}{a_{n+m-1}+\mu_1},~~\al_{n+m-1}(\ga_2) = \frac{1}{a_{n+m-1}+\mu_2},$$ 
with some $\mu_1, \mu_2 \in [0,1)$. 
Hence we have $$\left|
\log{\displaystyle\frac{\al_{n+m-1}(\ga_1)}{\al_{n+m-1}(\ga_2)}} \right| < \log{2}< 1 = 2^{0}.$$
\end{proof}

\medskip
\noindent
We now bound the influence of the difference on the 
$\log{\displaystyle\frac{1}{\al_i}}$ terms. 

\begin{lem}
\label{lem4}
Let $\ga_1$ and $\ga_2$ be two numbers whose continued fraction 
expansions coincide in the first $n+m-1$ terms $[a_1, a_2, \ldots, a_{n+m-1}]$. 
Then for any $i <  n+m-1$ we have 
$$
\left|
\log{\displaystyle\frac{\log \displaystyle\frac{1}{\al_{i}(\ga_1)}}{\log 
\displaystyle\frac{1}{\al_{i}(\ga_2)}}} \right|
< 2^{i-(n+m)+2}.
$$
In particular, this applies with $\ga_1 = \be^N$ and $\ga_2 = \be^1$. 
\end{lem}

\begin{proof}Assume that $\al_{i}(\ga_1)\le \al_{i}(\ga_2)$, the reverse 
case is 
done in the same way. In this case we need to prove 
$$
\displaystyle\frac{\log \displaystyle\frac{1}{\al_{i}(\ga_1)}}{\log
\displaystyle\frac{1}{\al_{i}(\ga_2)}} < e^{2^{i-(n+m)+2}}.
$$
Denote $c = \al_{i+1}(\ga_1)$ and $d =\al_{i+1} ( \ga_2)$. 
Then we have $\al_{i}(\ga_1) = \displaystyle\frac{1}{k+c}$ and 
$\al_{i}(\ga_2) = \displaystyle\frac{1}{k+d}$ for some integer $k \ge 1$. 
Hence $0 < d \le c < 1$. We have 
$$
\displaystyle\frac{\log \displaystyle\frac{1}{\al_{i}(\ga_1)}}{\log
\displaystyle\frac{1}{\al_{i}(\ga_2)}} = \displaystyle\frac{\log(k+c)}{\log(k+d)}.
$$
By Lemma \ref{lem2} we know that $\displaystyle\frac{c}{d} < e^{2^{i-(n+m)+2}}$, 
hence it suffices to show that $\displaystyle\frac{\log(k+c)}{\log(k+d)} \le 
\displaystyle\frac{c}{d}$. This is equivalent to $\displaystyle\frac{\log (k+c)}{c} \le 
\displaystyle\frac{\log (k+d)}{d}$. Consider the function $f(x) = \displaystyle\frac{\log(k+x)}{x}$
on the interval $(0,1)$. The reader can readily verify that 
$f'(x)<0$ for $x\in (0,1)$ so that $f$ is decreasing on this interval, and hence 
$f(c) \le f(d)$, which completes the proof.
\end{proof}

\noindent
We are now ready to bound the influence of changes in $N$ on the 
value of $\Phi^{-}$. 

\begin{lem}
\label{lem5}
For any $\om$ of the form as in Lemma \ref{smlchg} and for any $\ve >0$, 
there is an $m_0 >0$ such that for any $N$ and any $m \ge m_0$,
$$
 | \Phi^{-} (\be^{N}) -
\Phi^{-} ( \be^1) | < \displaystyle\frac{\ve}{4}.
$$
\end{lem}

\begin{proof}
The $\sum$ in the expression for $\Phi (\be^1)$ converges, 
hence there is an $m_1 >1$ such that the tail of the sum 
$\sum_{i\ge n+m_1} \al_1 \al_2 \ldots \al_{i-1} \log \displaystyle\frac{1}{\al_i} < 
\displaystyle\frac{\ve}{40}$. We will show how to choose $m_0 > m_1$ to
satisfy the conclusion of the lemma. 

We bound the influence of the change from $\be^1$ to 
$\be^N$ using Lemmas \ref{lem2} and \ref{lem4}. The influence on 
each of the ``head elements" ($i< n+m_1$) is bounded by 
$$
\left| \log{\displaystyle\frac{\al_1 (\be^1) \ldots \al_{i-1}(\be^1) \log 
\displaystyle\frac{1}{\al_i(\be^1)}}{\al_1 (\be^N) \ldots \al_{i-1}(\be^N) \log
\displaystyle\frac{1}{\al_i(\be^N)}}} 
\right| < \sum_{j=1}^{i-1} 2^{j-(n+m)+1} + 2^{i-(n+m)+2} < 2^{i -(n+m)+3}<
2^{m_1+3-m}.
$$
By making $m$ sufficiently large (i.e. by choosing a sufficiently large
$m_0$) we can ensure that 
$$
1- \displaystyle\frac{\ve}{40 \Phi(\be^1)} < \displaystyle\frac{\al_1 (\be^N) \ldots 
\al_{i-1}(\be^N) \log
\displaystyle\frac{1}{\al_i(\be^N)}}{\al_1 (\be^1) \ldots \al_{i-1}(\be^1) \log
\displaystyle\frac{1}{\al_i(\be^1)}} < 1 + \displaystyle\frac{\ve}{40 \Phi(\be^1)},
$$
hence
$$
\left| \al_1 (\be^N) \ldots
\al_{i-1}(\be^N) \log
\displaystyle\frac{1}{\al_i(\be^N)} - \al_1 (\be^1) \ldots
\al_{i-1}(\be^1) \log
\displaystyle\frac{1}{\al_i(\be^1)} \right| < 
$$
$$\displaystyle\frac{\ve}{40 \Phi(\be^1)}\al_1 (\be^1)
\ldots
\al_{i-1}(\be^1) \log
\displaystyle\frac{1}{\al_i(\be^1)}.
$$

Adding the inequality for $i=1,2, \ldots, n+m_1-1$ we obtain 
$$
\left| \sum_{i=1}^{n+m_1-1}  \al_1 (\be^N) \ldots
\al_{i-1}(\be^N) \log
\displaystyle\frac{1}{\al_i(\be^N)} - \sum_{i=1}^{n+m_1-1}  \al_1 (\be^1) \ldots
\al_{i-1}(\be^1) \log
\displaystyle\frac{1}{\al_i(\be^1)} \right| < 
$$
$$
\displaystyle\frac{\ve}{40 \Phi(\be^1)} \sum_{i=1}^{n+m_1-1}  \al_1 (\be^1) \ldots
\al_{i-1}(\be^1) \log
\displaystyle\frac{1}{\al_i(\be^1)} < \displaystyle\frac{\ve}{40 
\Phi(\be^1)} \Phi(\be^1) = 
\displaystyle\frac{\ve}{40}.
$$
Hence the influence on the ``head" of $\Phi^-$ is bounded 
by $\displaystyle\frac{\ve}{40}$. 

To bound the influence on the ``tail" we consider three kinds of 
terms \\$\al_1 (\be^N) \ldots
\al_{i-1}(\be^N) \log
\displaystyle\frac{1}{\al_i(\be^N)}$: those for which $n+m_1 \leq i \leq n+m-2$, $i=m+n-1$ and
$i \ge m+n+1$ (recall that $i=n+m$ is not in $\Phi^-$).

\medskip

\noindent
$\bullet$
{For $n+m_1 \leq i \leq n+m-2$. By Lemmas \ref{lem2} and \ref{lem4}:} 
$$
\left| \log{\displaystyle\frac{\al_1 (\be^1) \ldots \al_{i-1}(\be^1) \log
\displaystyle\frac{1}{\al_i(\be^1)}}{\al_1 (\be^N) \ldots \al_{i-1}(\be^N) \log
\displaystyle\frac{1}{\al_i(\be^N)}}}
\right| < \sum_{j=1}^{i-1} 2^{j-(n+m)+1} + 2^{i-(n+m)+2} < 2^{i -(n+m)+3}\le
2.
$$
Hence in this case each term can increase  by a factor of at most $e^2$. 

\medskip

\noindent
$\bullet$
{For $i=n+m-1$} Note that the change decreases 
$\log{\displaystyle\frac{1}{\al_{n+m-1}}}$
so that $\log{\displaystyle\frac{1}{\al_{n+m-1}(\be^N)}} \le 
\log{\displaystyle\frac{1}{\al_{n+m-1}(\be^1)}}$, hence we have 
$$
\log{\displaystyle\frac{\al_1 (\be^N) \ldots \al_{i-1}(\be^N) \log
\displaystyle\frac{1}{\al_i(\be^N)}}{\al_1 (\be^1) \ldots \al_{i-1}(\be^1) \log
\displaystyle\frac{1}{\al_i(\be^1)}}} \le 
\log{\displaystyle\frac{\al_1 (\be^N) \ldots \al_{i-1}(\be^N)}{\al_1 (\be^1) \ldots 
\al_{i-1}(\be^1)}} < 
$$
$$
\sum_{j=1}^{n+m-2} 2^{j-(n+m)+1} < 1. 
$$
Hence this term could increase by a factor of at most ${e}$. 

\medskip

\noindent
$\bullet$
{For $ i \geq n+m+1$:}
Note that $\al_j$ for $j>n+m$ are not affected by the change, 
and the change decreases $\al_{n+m}$, so that $\al_{n+m} (\be^N) \le 
\al_{n+m} (\be^1)$. 
Hence 
$$
\log{\displaystyle\frac{\al_1 (\be^N) \ldots \al_{i-1}(\be^N) \log
\displaystyle\frac{1}{\al_i(\be^N)}}{\al_1 (\be^1) \ldots \al_{i-1}(\be^1) \log
\displaystyle\frac{1}{\al_i(\be^1)}}} = 
\log{\displaystyle\frac{\al_1 (\be^N) \ldots \al_{n+m}(\be^N)}{\al_1 (\be^1) \ldots 
\al_{n+m}(\be^1)}} \le
$$
$$
\log{\displaystyle\frac{\al_1 (\be^N) \ldots \al_{n+m-1}(\be^N)}{\al_1 (\be^1) \ldots
\al_{n+m-1}(\be^1)}} < \sum_{j=1}^{n+m-1} 2^{j-(n+m)+1} < 2
$$
So in this case each term could increase by a factor of at most $e^2$. 

We see that after the change each term of the tail could increase 
by a factor of $e$ at most. The value of the tail remains 
positive in the interval $\left(0, \displaystyle\frac{e^2 \ve}{40}\right]$, hence the 
change in the tail is bounded  by $\displaystyle\frac{e^2 \ve}{40} < \displaystyle\frac{ 9 \ve}{40}$. 

\smallskip
\noindent
So the total change in $\Phi^-$ is bounded by
$$
\mbox{change in the ``head"}~~ +~~
\mbox{change in the ``tail"} ~~ < \displaystyle\frac{\ve}{40} + \displaystyle\frac{9 \ve}{40} = 
\displaystyle\frac{\ve}{4}.
$$
\end{proof}

\noindent
The following Lemma follows immediately from Lemma \ref{lem5}. 

\begin{lem}
\label{lem6}
For any $\ve$ and for the same $m_0 (\ve)$ as in Lemma \ref{lem5}, 
for any $m \ge m_0$ and $N$, 
$$
 | \Phi^{-} (\be^{N}) -
\Phi^{-} ( \be^{N+1}) | < \displaystyle\frac{\ve}{2}.
$$ 
\end{lem}

\begin{proof}
We have 
$$
| \Phi^{-} (\be^{N}) -
\Phi^{-} ( \be^{N+1}) | \le | \Phi^{-} (\be^{N}) -
\Phi^{-} ( \be^1) | + | \Phi^{-} (\be^{1}) -
\Phi^{-} ( \be^{N+1}) | < \displaystyle\frac{\ve}{4} + \displaystyle\frac{\ve}{4} = \displaystyle\frac{\ve}{2}.
$$
\end{proof}

\noindent
We will now have to take a closer look at the term 
$\al_1 \ldots \al_{n+m-1} \log \displaystyle\frac{1}{\al_{m+m}} = 
\Phi(\om) - \Phi^- (\om)$. 
We will need the following simple statement. 

\begin{lem}
\label{lem8}
For any $k>1$, $\al_{k-1} \al_{k} < \displaystyle\frac{1}{2}$.
\end{lem}
\begin{proof}
There is an integer $l \geq 1$ such that 
$$
\al_{k-1} \al_{k} = \displaystyle\frac{1}{l+\al_k} \al_{k} <
\displaystyle\frac{1}{\al_k + \al_k} \al_k = \displaystyle\frac{1}{2}.
$$
\end{proof}

\noindent
Denote $\Phi^1 (\al) = \al_1 \ldots \al_{n+m-1} \log \displaystyle\frac{1}{\al_{n+m}} =
\Phi(\om) - \Phi^- (\om)$, we are now ready to prove the following. 

\begin{lem}
\label{lem9}
For sufficiently large $m$, for any $N$, 
$$
\Phi^1 (\be^{N+1}) - \Phi^1 (\be^N) < \displaystyle\frac{\ve}{2}.
$$
\end{lem}

\begin{proof}
According to \lemref{lem7} we have 
$$
\left|
\log{\displaystyle\frac{\al_{1}(\be^{N+1})\ldots 
\al_{n+m-1}(\be^{N+1})}{\al_{1}(\be^{N}) \ldots \al_{n+m-1}(\be^{N})}} 
\right| < \sum_{i=1}^{n+m-1} 2^{i-(n+m)}/N < \displaystyle\frac{1}{N}. 
$$
Hence $\al_{1}(\be^{N+1})\ldots
\al_{n+m-1}(\be^{N+1}) < \al_{1}(\be^{N}) \ldots \al_{n+m-1}(\be^{N}) 
e^{1/N}$, and 
$$
\Phi^1 (\be^{N+1}) < \Phi^1 (\be^N) e^{1/N}
\displaystyle\frac{\log{\displaystyle\frac{1}{\al_{n+m}(\be^{N+1})}}}
     {\log{\displaystyle\frac{1}{\al_{n+m}(\be^{N})}}} = 
\Phi^1 (\be^N) e^{1/N}
\displaystyle\frac{\log (N+1+1/\phi)}{\log(N+1/\phi)}.
$$
Hence 
$$
\Phi^1 (\be^{N+1}) - \Phi^1 (\be^N) < \Phi^1 (\be^N) \left(e^{1/N}
\displaystyle\frac{\log (N+1+1/\phi)}{\log(N+1/\phi)}-1\right) <
$$
$$
\Phi^1 (\be^N)\left(\left(1+\displaystyle\frac{e}{N}\right) \displaystyle\frac{\log 
(N+1+1/\phi)}{\log(N+1/\phi)}-1\right).
$$
We make the following calculations.
Denote $x = \displaystyle\frac{\log
(N+1+1/\phi)}{\log(N+1/\phi)}$, then $(N+1/\phi)^x = N+1+1/\phi$, and
$$(N+1/\phi)^{x-1} =  \displaystyle\frac{N+1+1/\phi}{N+1/\phi}= 1 + \displaystyle\frac{1}{N+1/\phi}< e^{\displaystyle\frac{1}{N+1/\phi}}.$$ 
$N+1/\phi > e^{1/3}$, and so $x-1 < \displaystyle\frac{3}{N+1/\phi} < \displaystyle\frac{3}{N}$,
thus $x < 1 + \displaystyle\frac{3}{N}$. 

\medskip
\noindent
By \lemref{lem8} we have 
$$
\Phi^1 (\be^N) = \al_1 (\be^N) \ldots \al_{n+m-1} (\be^N) \log 
\displaystyle\frac{1}{\al_{n+m} (\be^N)} < \left( \displaystyle\frac{1}{2} \right)^{(n+m-2)/2} 
\log ( N + 1/\phi). 
$$
Thus 
$$
\Phi^1 (\be^{N+1}) - \Phi^1 (\be^N) <
\Phi^1 (\be^N)\left(\left(1+\displaystyle\frac{e}{N}\right) \displaystyle\frac{\log
(N+1+1/\phi)}{\log(N+1/\phi)}-1\right) <
$$
$$
\left( \displaystyle\frac{1}{2} \right)^{(n+m-2)/2}
\log ( N + 1/\phi) \left( (1+e/N)(1+3/N) -1 \right) < 
\left( \displaystyle\frac{1}{2} \right)^{(n+m-2)/2} \log ( N + 1/\phi)
\displaystyle\frac{14}{N}.
$$
Since $\displaystyle\frac{14}{N} \in o(1/\log ( N + 1/\phi))$, this expression 
can be always made less than $\displaystyle\frac{\ve}{2}$ by choosing $m$ 
large enough. 
\end{proof}

\noindent
Lemmas \ref{lem6} and \ref{lem9} yield the following

\begin{lem}
\label{lem10}
For sufficiently large $m$, for any $N$, 
$$
\Phi (\be^{N+1}) - \Phi (\be^N) < \ve.
$$
\end{lem}
\begin{proof}
We use Lemmas \ref{lem6} and \ref{lem9}. For sufficiently large $m$, 
$$\Phi (\be^{N+1}) - \Phi (\be^N) \leq 
\Phi^- (\be^{N+1}) - \Phi^- (\be^N) + \Phi^1 (\be^{N+1}) - \Phi^1 (\be^N) 
< \displaystyle\frac{\ve}{2} + \displaystyle\frac{\ve}{2} = \ve. 
$$
\end{proof}

\noindent
To complete the proof of Lemma \ref{smlchg} we will need the 
following statement. 
\begin{lem}
\label{lem11}
$$
\lim_{N \rightarrow \infty} \Phi (\be^N) = \infty.$$
\end{lem}
\begin{proof}
We will prove that $\lim_{N \rightarrow \infty} \Phi^1 (\be^N) = \infty$,
this suffices, since $\Phi^1 (\be^N) < \Phi (\be^N)$.
By \lemref{lem2},
$$
\left|
\log{\displaystyle\frac{\al_{1}(\be^{N})\ldots 
\al_{n+m-1}(\be^{N})}{\al_{1}(\be^{1})\ldots
\al_{n+m-1}(\be^{1})}} \right| < \sum_{i=1}^{n+m-1} 2^{i-(n+m)+1} < 2,
$$ 
hence 
$$
\al_{1}(\be^{N})\ldots
\al_{n+m-1}(\be^{N}) > \displaystyle\frac{1}{e^2} \cdot \al_{1}(\be^{1})\ldots
\al_{n+m-1}(\be^{1})
$$
and
$$
\Phi^1 (\be^N) > \displaystyle\frac{1}{e^2} \cdot \displaystyle\frac{\log(N+1/\phi)}{\log(1+1/\phi)} \Phi^1 
(\be^1).
$$
The latter expression obviously goes to $\infty$ as $N \rightarrow 
\infty$. 
\end{proof}

\noindent
We are now ready to prove \lemref{smlchg}. 

\begin{proof} (of \lemref{smlchg}). Choose $m$ large enough for Lemma 
\ref{lem10} to hold. Increase $N$ by one at a time starting with $N=1$. 
We know that $\Phi(\be^1) = \Phi(\al) < \Phi(\al) + \ve$, and 
by  \lemref{lem11}, there exists an $M$ with $\Phi(\be^M) > \Phi(\al) + 
\ve$. Let $N$ be the smallest such $M$.  Then $\Phi(\be^{N-1}) \le
\Phi(\al) +\ve$, and by  \lemref{lem10}
$$
\Phi(\be^N) < \Phi(\be^{N-1}) + \ve \le \Phi(\al) + 2 \ve.
$$
Hence 
$$
\Phi(\al) +\ve < \Phi(\be^N) <  \Phi(\al) + 2 \ve.
$$
Choosing $\be = \be^N$ completes the proof. 

\end{proof}

\noindent
We will now prove Lemma \ref{notdeclem}.

\begin{proof} (of \lemref{notdeclem}).   
The $\sum$ in the expression for $\Phi (\om)$ converges, 
hence there is an $m_1 >1$ such that the tail of the sum 
$\sum_{i\ge n+m_1} \al_1 \al_2 \ldots \al_{i-1} \log \displaystyle\frac{1}{\al_i} < 
\displaystyle\frac{\ve}{2}$. We will show how to choose $m_0 > m_1$ to
satisfy the conclusion of the lemma. 

By Lemmas \ref{lem2} and \ref{lem4}, for any $\be^T$ and any $i\le n+m_1$ we have
$$
\left| \log{\displaystyle\frac{\al_1 (\be^T) \ldots \al_{i-1}(\be^T) \log
\displaystyle\frac{1}{\al_i(\be^T)}}{\al_1 (\om) \ldots \al_{i-1}(\om) \log
\displaystyle\frac{1}{\al_i(\om)}}}
\right| < \sum_{j=1}^{i-1} 2^{j-(n+m)+1} + 2^{i-(n+m)+2}$$ $$ < 2^{i -(n+m)+3}\le
2^{n+m_1 -(n+m_0)+3}=2^{m_1-m_0+3}.
$$
We can choose $m_0$ sufficiently large so that $\displaystyle{e^{-2^{m_1-m_0+3}}}>1-\displaystyle\frac{\ve}{2\Phi(\om)}$. So 
that 
$$
\al_1 (\be^T) \ldots \al_{i-1}(\be^T) \log
\displaystyle\frac{1}{\al_i(\be^T)} > \left( 1-\displaystyle\frac{\ve}{2\Phi(\om)}\right)
\al_1 (\om) \ldots \al_{i-1}(\om) \log
\displaystyle\frac{1}{\al_i(\om)},
$$
for $i\le n+m_1$. 

Now, for any $\be^T$ we have
$$
\Phi(\be^T)   \ge
\sum_{i=1}^{n+m_1-1} \al_1(\be^T) \al_2(\be^T) \ldots \al_{i-1}(\be^T) \log \displaystyle\frac{1}{\al_i(\be^T)}>
$$
$$\sum_{i=1}^{n+m_1-1}
\left( 1-\displaystyle\frac{\ve}{2\Phi(\om)}\right)
\al_1 (\om) \ldots \al_{i-1}(\om) \log
\displaystyle\frac{1}{\al_i(\om)} = 
$$
$$
\left( 1-\displaystyle\frac{\ve}{2\Phi(\om)}\right) \left( \Phi(\om) - \sum_{i=n+m_1}^{\infty}
\al_1 (\om) \ldots \al_{i-1}(\om) \log
\displaystyle\frac{1}{\al_i(\om)} \right) >
$$
$$
\left( 1-\displaystyle\frac{\ve}{2\Phi(\om)}\right) \left( \Phi(\om) - \displaystyle\frac{\ve}{2} \right) >
\Phi(\om)-\ve.
$$
\end{proof}

\noindent
We will also need the 
following lemma in the proof of the Main Theorem.

\begin{lem}
\label{tailto1}
Let $\om = [a_1, a_2, a_3, \ldots]$ and let $\ve>0$ be given. 
Then there is an $N = N(\ve)$ such that for any $n \ge N$ 
we have $\Phi(\om_n)<\Phi(\om)+\ve$, where 
$\om_n = [a_1, a_2, \ldots, a_{n}, 1, 1, 1, \ldots]$. 
\end{lem}

\noindent
The proof is not hard and is similar to the proof of Lemma \ref{lem5}.
We present the main steps in the proof. 

\noindent
$\bullet$
There is an $m_0$ such that the sum of the tail elements of 
$\Phi(\om)$
is small:
$$
\sum_{i=m_0}^{\infty} \al_1 (\om) \ldots \al_{i-1} (\om) \log 
\frac{1}{\al_i(\om)} < \frac{\ve}{4 \cdot e^2}.
$$

\noindent
$\bullet$
Similarly to Lemma \ref{lem5}, we can use Lemmas \ref{lem2} and 
\ref{lem4} to show that for sufficiently large  $m_1> m_0$, $n > m_1$ 
implies 
that 
$$
{\displaystyle\frac{\al_1 (\om_n) \ldots \al_{i-1}(\om_n) \log 
\displaystyle\frac{1}{\al_i(\om_n)}}{\al_1 (\om) \ldots \al_{i-1}(\om) 
\log\displaystyle\frac{1}{\al_i(\om)}}} < 1+ \frac{\ve}{4 \Phi(\om)}
$$
for all $i < m_0$. 

\noindent
$\bullet$
Again by Lemmas \ref{lem2}, \ref{lem4} and \ref{lem8} to show that for any 
$i$ (with 
a special consideration 
to the case $i=n$), 
$$
{\al_1 (\om_n) \ldots \al_{i-1}(\om_n) \log 
\displaystyle\frac{1}{\al_i(\om_n)}}
< \max \left( e^2 \cdot \al_1 (\om) \ldots \al_{i-1}(\om) 
\log\displaystyle\frac{1}{\al_i(\om)},2^{2-i/2}\right).
$$

\noindent
Adding these up we get for $n>m_1$:
$$
\Phi(\om_n) < \left( 1+ \frac{\ve}{4 \Phi(\om)}\right) \sum_{i=1}^{m_0-1}
 \al_1 (\om) \ldots \al_{i-1} (\om) \log \frac{1}{\al_i(\om)} +
$$
$$ 
\sum_{i=m_0}^{\infty}  \max \left( e^2 \cdot \al_1 (\om) \ldots 
\al_{i-1}(\om)
\log\displaystyle\frac{1}{\al_i(\om)},2^{2-i/2}\right) \le 
$$
$$
 \left( 1+ \frac{\ve}{4 \Phi(\om)}\right)  \Phi(\om) + 
\sum_{i=m_0}^{\infty} e^2 \cdot \al_1 (\om) \ldots
\al_{i-1}(\om)
\log\displaystyle\frac{1}{\al_i(\om)}
+ \sum_{i=m_0}^{\infty}
2^{2-i/2} \le
$$
$$
\Phi(\om) + \frac{\ve}{4} + \frac{e^2 \cdot \ve}{4 \cdot e^2} + 
2^{4-m_0/2}=
\Phi(\om) + \frac{\ve}{2} + 
2^{4-m_0/2}.
$$
We complete the proof by choosing $m_0$ large enough so that $2^{4-m_0/2} 
<\ve/2$. 

\section{Proof of \thmref{no partition}}

\comm{

\noindent
We are now ready to prove the Main Theorem. 

\medskip

\noindent
\begin{thm}
\label{almostmain}
There exists a parameter value $\theta\in\RR/\ZZ$ such that the Julia set of the
quadratic polynomial $f_\theta (z)=z^2+e^{2\pi i \theta}z$ is not computable.
\end{thm}

\medskip

\noindent
We note the the Main Theorem follows from
\thmref{almostmain} in an elementary way. Indeed,
let $\theta \in \RR/\ZZ$ be as promised by \thmref{almostmain}. Then given 
an oracle access to $\theta$ it is  easy to compute $\lambda = e^{2 \pi i \theta}$, 
and vice versa. This implies that we could simulate the oracle for $\lambda$ given an oracle for 
$\theta$, and thus it is impossible to compute $J(z^2 + \lambda z)$ with an oracle 
access to $\lambda$.

\noindent
Setting $c=\lambda^2/4-\lambda/2$ and $T_c(z)=z-2\lambda$, we have 
$$J(z^2 + \lambda z) = T_c(J(z^2 + c)).$$ Hence being able to compute 
$J(z^2 + c)$ with an oracle access to $c$ would allow us to compute 
$J(z^2 +\lambda z)$ with an oracle access to $\lambda $.

\medskip

\noindent
{\bf Proof of Theorem \ref{almostmain}.}
There are countably many oracle Turing Machines,
hence we can enumerate them as $M_1^{\phi}, M_2^{\phi}, \ldots$. 
We prove the theorem by constructing a 
number $c$ that ``fools" all the oracle Turing Machines $M_1^{\phi}, M_2^{\phi}, \ldots$
which are trying to compute $J_c = J(f_c)$. 
We proceed iteratively, so on step $i$ we ``fool" the first $i$ machines $M_k^{\phi}$. 
We then take the process to the limit, obtaining a number on which none 
of the machines works. 

}

Recall that $r(\theta)$ denotes the conformal radius of the Siegel disk associated
with the polynomial $P_\theta (z) = z^2 + e^{2 \pi i \theta} z$, or zero,
if $\theta$ is not a Brjuno number.

We will argue by way of contradiction, and assume that 
there exists a countable union of sets
$$\cup_{i=1}^\infty S_i=\RR/\ZZ$$
such that the function $r(\theta)$ is continuous on each $S_i$.

\begin{mlem} 
\label{Main:Lemma}
There exist
\begin{itemize}
\item
a sequence of initial segments $I_i = [a_1, a_2, \ldots, a_{N_i}]$, and
\item 
a sequence of nested intervals 
$$[l_0, r_0] \supset [l_1, r_1] \supset [l_2, r_2] \supset \ldots,$$
\end{itemize}
such that the following properties are maintained:
\begin{enumerate}
\item
\label{cond1}
whenever $i>j$ we have
$$I_i=[I_j,a_{N_j+1},a_{N_j+2},\ldots,a_{N_i}];$$
\item 
\label{cond2}
$r_i = r(\ga_i)$, where $\ga_i = [I_i, 1,1, \ldots]$;
\item 
\label{cond3}
for each $i\ge 1$ and for every
 $\be=[I_i, t_{N_i+1}, t_{N_i+2}, \ldots]$ with $r(\be)\in [l_i, r_i]$ 
we have $$\be\notin S_i;$$
\item 
\label{cond4}
denote $\ell_i = r_i - l_i$. Then 
$$\ell_i>0\text{ and }\ell_{i}\le\ell_{i-1}/2\text{ for all }i\geq 1;$$
\item
\label{cond5}
for any $\be=[I_i, t_{N_i+1}, t_{N_i+2}, \ldots]$, $i\ge 1$, we have
$$\Phi(\be)>\Phi(\ga_{i-1})-2^{-(i-1)}.$$
\end{enumerate}
\end{mlem}

\noindent
\begin{proof}[Proof of the Main Lemma.] We prove the Main Lemma by induction on $i$.
For the basis of induction, set $I_0 = [1]$, $r_0= r(\ga_0)$ and $l_0 = r_0/2$, where
$\ga_0 = [1,1,1, \ldots]$. Then for $i=0$ conditions (1)-(5) trivially hold.

\subsection*{The induction step.} We now have the conditions satisfied
 for some $i\geq 0$ and would like to extend them to $i+1$. 

Set $S\equiv S_{i+1}$ and
let $R$ be
the set of all possible values of the conformal radius $r(\theta)$ for $\theta\in S$. There
are two possibilities:

\medskip
\noindent
{\bf Case 1: }
There exist $\ve_0>0$ and $m_0\in\NN$ such that for every $\be\in S$ of the form 
$$\be=[I_i,\underbrace{1,1,\ldots,1}_{m_0},\ldots]\text{ we have }|r_i-r(\be)|>\ve_0.$$
In this case, select $0<\ve\leq\min(\ve_0,\ell_i/2)$.
Set
 $$I_{i+1}=[I_i,\underbrace{1,1,\ldots,1}_{m_0}],\;l_{i+1} = r_i - \ve,\text{ and }
r_{i+1} = r_i.$$  
$\ga_{i+1} = [I_i,1, 1, \ldots ] = \ga_i$. We have
$r(\ga_{i+1}) = r(\ga_i) = r_i = r_{i+1}$, hence  conditions 
\eref{cond1},\eref{cond2} and \eref{cond4} are satisfied. 

Suppose $\be = [I_{i+1}, t_{N_{i+1}+1}, t_{N_{i+1}+2}, \ldots]$ with $r(\be) \in 
[l_{i+1}, r_{i+1}]$. 
Then $\be\notin S$ and 
\eref{cond3} is satisfied. By Lemma \ref{notdeclem}, we can choose
$m_0$ sufficiently large in $I_{i+1}$, so that for any $\be$ beginning 
with 
$I_{i+1}$, we have $\Phi(\be)>\Phi(\ga_i)-2^{-i}$ thus satisfying \eref{cond5}. 

\medskip
\noindent
The complementary case is the main part of the argument:

\noindent
{\bf Case 2.} 
For every $\ve>0$ and $m\in\NN$ we can find $\beta\in S$ starting with $I_i$ followed by $m$ ones
so that
\begin{equation}
\label{be1}
r_i - \ve< r(\be) \le r_i
\end{equation}

Choose an $\ve$ such that $r_i - 3 \ve > l_i > r_i - 4 \ve$. Denote 
$$\ve_0 = \min \left( \displaystyle\frac{\log(r_i-\ve)-\log(r_i- 2 \ve)}{8},
\frac{\log(r_i-2 \ve)-\log(r_i- 3 \ve)}{8} \right)>0.$$

\noindent
Theorem \ref{phi-cont} of Buff and Ch{\'e}ritat 
says that the function 
$$\upsilon:\theta \mapsto \Phi(\theta)+\log r(\theta)$$
 continuously extends to $\RR/\ZZ$. 
Due to compactness of $\RR/\ZZ$, this function is uniformly continuous, 
and there exists a $\de_0>0$ such that if $|x-y|< \de_0$ then $|\upsilon(x)-\upsilon(y)|<\ve_0$. 

\medskip

\noindent
We choose $m$ large enough, so that for any 
$\zeta = [I_i, \underbrace{1,1, \ldots, 1}_{m}, \ldots]$ we have
$| \ga_i - \zeta |< \de_0$ and so that Lemma \ref{notdeclem} holds for 
 $m_0=m$ with $I=I_i$ and $\ve=2^{-i}$. 
Write $$\be = [I_i, \underbrace{1,1, \ldots, 1}_{m},  t_{N_i+m+1}, 
t_{N_i+m+2}, \ldots]\in S.$$ 

\noindent
By assumption, the conformal radius $r(\bullet)$ is continuous on $S$. 
Hence there is a
$\de>0$ such that $$|  r(x) -  r(\be)| < \ve\text{ whenever } 
|x - \be| < \de\text{ and }x \in S.$$

\noindent
By Lemma \ref{tailto1}, there is 
an $N$ such that for any $n \ge N$,  $$\be_n = [I_i, \underbrace{1,1, 
\ldots, 1}_{m}, t_{N_i+m+1}, \ldots, t_{N_i+m+n}, 1, 1, \ldots]$$
satisfies $$\Phi(\be_n) < \Phi(\be) + \ve_0.$$  
We can choose $n \ge N$ large 
enough so that for any $x$ whose contnued fraction expansion has the initial segment 
$$I_i^0 = [I_i, \underbrace{1,1, \ldots, 1}_{m}, t_{N_i+m+1}, t_{N_i+m+2}, 
\ldots , t_{N_i+m+n}],$$ we have
$$|x - \be| < \min(\de_0/2,\de).$$

\smallskip

Start with $\om_0 = \be_n = [I_i^0, 1, 1, \ldots]$. We have 
$| \om_0 - \be | < \de_0$, and hence $|\upsilon(\om_0)-\upsilon(\be)|< \ve_0$. So

$$
\log r(\om_0) = \upsilon(\om_0) - \Phi(\om_0 ) > \upsilon(\be) - \ve_0 - \Phi(\be)- \ve_0= 
\log r(\be) - 2 \ve_0.
$$

\noindent
 By Lemma \ref{smlchg} we can 
extend $I_i^0$ to a longer initial segment $I_i^1$ so that 
setting $\om_1= [I_i^1, 1, 1, \ldots]$ we have  $$\Phi(\om_0)+2 \ve_0< 
\Phi(\om_1) < \Phi(\om_0)+4 \ve_0.$$

We have
 $|\om_0-\be|<\de_0/2$ and $|\om_1-\be|<\de_0/2$, so 
$|\om_0 - \om_1|<\de_0$, and $|\upsilon(\om_0)-\upsilon(\om_1)|<\ve_0$. Hence 
$$
\log (r(\om_1)) = \upsilon(\om_1) - \Phi(\om_1) > \upsilon(\om_0) - \ve_0 - \Phi(\om_0)
- 4 \ve_0 = \log (r(\om_0)) - 5 \ve_0,
$$
and
$$
\log (r(\om_1)) = \upsilon(\om_1) - \Phi(\om_1) < \upsilon(\om_0) + \ve_0 - \Phi(\om_0)
- 2 \ve_0 = \log (r(\om_0)) -  \ve_0.
$$
Hence $$\log (r(\om_0)) -   5 \ve_0 < \log (r(\om_1)) < \log (r(\om_0)) -  \ve_0.$$

\smallskip

\noindent
In the same fashion, we can extend $I_i^1$ to $I_i^2$, and obtain 
$\om_2= [I_i^2, 1, 1, \ldots]$ so that 
$$\log (r(\om_1)) -   5 \ve_0 < \log (r(\om_2)) < \log (r(\om_1)) -  \ve_0.$$

\medskip
\noindent
Recall that $$\log (r(\om_0)) > \log (r(\be))- 2 \ve_0 > \log (r_i -  \ve) - 2 \ve_0
\ge  \log (r_i -  2 \ve) + 6 \ve_0.$$ Hence,
after finitely many steps, we will obtain $I_i^k$ and $\om_k= [I_i^k, 1, 1, \ldots]$
such that 
$$
\log ( r_i - 3 \ve ) + \ve_0 < \log ( r(\om_k)) < \log (r_i - 3 \ve ) + 6 \ve_0 < 
\log(r_i - 2 \ve). 
$$
Choose $I_{i+1} = I_i^k$, $\ga_{i+1} = \om_k$,  $l_{i+1}=l_i$, and $r_{i+1} = r(\om_k)$. We have 
$$l_{i+1} < r_i - 3\ve < r_{i+1} < r_i - 2 \ve.$$
 Condition \eref{cond1} and \eref{cond2} are satisfied by definition.
Condition \eref{cond4} is satisfied because $$\ell_i =r_i-l_i< 4\ve\text{ and }\ell_{i+1} <
\ell_i - 2\ve < \ell_i/2.$$
Condition \eref{cond3} is satisfied, because 
$x = [I_i^0, \ldots]\in S$ only when  $r(x) \in 
[r(\be)-\ve,r(\be)+\ve] \subset$ $[r_i-2 \ve, r_i + \ve]$, 
and by our construction $[l_{i+1}, r_{i+1}]$ is disjoint from  $[r_i-2 
\ve, r_i + \ve]$.
Finally, condition \eref{cond5} is satisfied because Lemma \ref{notdeclem} holds with $\ve=2^{-i}$ with 
the $m$ chosen in the initial segment $I_i^0$. We have thus completed the proof of the Main Lemma \ref{Main:Lemma}.
 
\end{proof}

\subsection*{Passing to the Limit}
The completion of the proof of \thmref{no partition} relies  on the 
following lemma:

\begin{lem}
\label{tothelimit}
Denote $\ga = \lim_{i \rightarrow \infty} \ga_i$. Then the following equalities hold:
$$
\Phi(\ga) = \lim_{i \rightarrow \infty}\Phi(\gamma_i)~~~~\mbox{$~$and$~$}~~~~
r(\ga) = \lim_{i \rightarrow \infty}r(\gamma_i).
$$
\end{lem}

\begin{proof}
By the construction, the limit $\ga = \lim \ga_i$ exists. We also 
know by condition \eref{cond4} of the Main Lemma \ref{Main:Lemma}
that the sequence $r(\ga_i)=r_i$ converges uniformly to 
some number $r$. By condition \eref{cond5} of the Main Lemma \ref{Main:Lemma}
the sequence $\Phi(\ga_i)-2^{-(i-1)}$ is 
 non-decreasing, and hence converges to a value $\psi$ ({\em a priori}
we could have $\psi=\infty$). The sequence $\Phi(\ga_i)$ must converge to $\psi$ as
well. 


By \propref{r doest drop}, we have $r(\ga) \ge r >0$.
On the other hand, by condition \eref{cond5} of the Main Lemma \ref{Main:Lemma}, 
we know that $\Phi(\ga)>\Phi(\ga_i)-2^{-(i-1)}$ for all $i$. 
Hence $\Phi(\ga) \ge \psi$. In particular $\psi<\infty$. 

From \cite{BC} we know that 
\beq
\label{BCCont}
\psi + \log r= \lim (\Phi(\ga_i) + \log r(\ga_i)) = \Phi(\ga) + \log r(\ga).
\eeq
Along with $r(\ga)\ge r$ and $\Phi(\ga)\ge \psi$ this yields 
 $\Phi (\ga)=\psi$, and $r(\ga)=r$, which completes
the proof. 
\end{proof}

\subsection*{Finalizing the argument.} Let $\ga$ be the limit from 
the previous lemma. We claim that $\ga\notin\cup S_i$.
 Indeed, for every $i$, the continued fraction $I_i$ is an initial segment of the continued
fraction expansion of $\ga$ by condition \eref{cond1} of the Main Lemma \ref{Main:Lemma}. 
By Lemma \ref{tothelimit} $$r(\ga) = \lim r(\ga_i)= \lim r_i \in [l_i,r_i].$$ 
Thus by condition \eref{cond3} of the Main Lemma \ref{Main:Lemma} we have $\ga\notin S_i$.
We have in this way arrived at a contradiction with $\cup S_i=\RR/\ZZ$, which completes the 
proof of \thmref{no partition}.

\section{Concluding remarks}
\subsection*{Connection with the work of Buff and Ch{\'e}ritat}
Let us outline here how the methods of \cite{BC2} can be applied to prove \thmref{no partition}
instead of the estimates of \secref{sec changes} (we note that a newer version of the 
same result exists \cite{ABC}, where the arguments we quote are simplified).
The main technical result of that paper is the following. Let $\alpha=[a_0,a_1,\ldots]$
be a Brjuno number, and as before denote $p_k/q_k$ the sequence of its continued fraction
approximants. Let $A>1$ and for each integer $n\geq 0$ set 
$$\alpha[n]=[a_0,a_1,\ldots,a_n,A^{q_n},1,1,1,\ldots].$$
Then for this  particular sequence of Brjuno approximants of $\alpha$,
$$\Phi(\alpha[n])\underset{n\to\infty}{\longrightarrow}\Phi(\alpha)+\log A,$$
and moreover,
$$\lim r(\alpha[n])= r(\alpha)/A.$$
The last equality can be used to construct the ``drops'' in the value of the conformal radius
of the Siegel disk needed to inductively avoid the sets $S_i$. In this way, one obtains a sequence
of Brjuno numbers $\theta_i\to\theta$ with conformal radii $r_i=r(\theta_i)>r_{i+1}$
such that $\lim r_i=r>0$, and $\theta_i$ is not in any of the $S_j$ up to $i$-th.

It remains to show that $r(\theta)=r$, as {\it a priori} only the inequality ``$\leq$'' is known. 
Buff and Ch{\'e}ritat demonstrate it in their context. The idea is, roughly speaking,
in showing that the boundary of $\Delta(\theta_i)$ is well approximated by a periodic cycle
of a high period. The perturbation $\theta_i\mapsto\theta_{i+1}$ is then chosen sufficiently
small so that the cycle does not move much.

As a final remark, let us point out:

\begin{rem}
Combining the methods of \cite{BC2} with our argument as outlined here one may strengthen
the Main Theorem by showing that there exists a non-computable Siegel Julia set for
which the boundary of the Siegel disk is {\it smooth.}
\end{rem}

\subsection*{Further progress} It is a natural question to ask whether the construction of non-computable
Julia sets carried out in this paper can be replaced with a different, perhaps, simpler approach.
Jointly with I.~Binder, we have demonstrated the following in \cite{BBY1}:

\begin{thm}[\cite{BBY1}]
\label{bby}
Let $R$ be a rational mapping of the Riemann sphere with no rotation domains (either Siegel disks
or Herman rings). Then its Julia set is computable by a TM $M^\phi$ with an oracle for the
coefficients of $R$. 

\end{thm}

\noindent
Moreover, it is shown in the same paper that the Julia set of a quadratic polynomial $J_c$ with
a periodic Siegel disk with conformal radius $r$ is computable by a TM with an oracle for $c$ if and only if 
$r$ itself is computable by some such machine. In retrospect, therefore, our approach finds the only
available class of  examples.

The size of the set of parameter values $\theta\in\RR/\ZZ$ for which $J(P_\theta)$ is uncomputable is rather meagre.
One can show combining the results of \cite{BBY1} with, for example, those of
Petersen and Zakeri \cite{PZ} that this set has Lebesgue measure zero;
and \thmref{bby} implies that its complement contains a dense $G_\delta$ subset of $\RR/\ZZ$.
It is natural to ask if, for example, its Hausdorff dimension is positive, and the answer to this question is 
not known to us. It is also interesting to ask if any values of $\theta$ in this
set are computable reals (as there are only countably many computable reals, and our procedure
clearly produces an uncountable set of $\theta$'s, most of them cannot be computable).
We again do not know the answer to this.

On the practical side of things, to our knowledge, one has not been able to produce informative pictures
of quadratic Julia sets with Cremer orbits, although by \thmref{bby} this is theoretically possible.
One potential explanation is that the computational complexity of these sets (the amount of time
it takes to decide whether to color a pixel of size $2^{-n}$ as a function of $n$) is very high.
This is indeed so for the na{\"\i}ve algorithms. In \cite{BBY2} jointly with I.~Binder we 
have constructed quadratic Julia sets whose computational complexity is arbitrarily high, but
again all with Siegel disks. 

A natural first step towards studying the complexity of Cremer Julia sets is to look at parabolics,
but the first author has recently demonstrated in \cite{Brv2} that having a parabolic orbit does not
qualitatively change the complexity of computing a Julia set. This opens an entertaining possibility
that some Cremer Julia sets have attainable computational complexity, and could be practically drawn by
a clever algorithm.


\begin{thebibliography}{*****}
\bibitem[Brv2]{Brv2} M. Braverman, Parabolic Julia sets are polynomial time computable. e-print math.DS/0505036 at Arxiv.org.
\end{thebibliography}
\end{document}